\documentclass[preprint]{elsarticle}
\usepackage[utf8]{inputenc}
\usepackage{amssymb}
\usepackage{amsmath}
\usepackage{amsthm}
\usepackage{hyperref}
 \usepackage{a4wide}
\usepackage{breqn}
\usepackage[table,xcdraw]{xcolor}
\usepackage{xcolor}
\usepackage{mathtools}
\usepackage{bm}
\usepackage{enumitem}
\usepackage{mathrsfs}
\usepackage{graphicx}
\usepackage{caption}
\usepackage{url}
\usepackage{euscript}
\usepackage{placeins}
\usepackage{mathdots}
\usepackage[english]{babel}
\usepackage{array}
\usepackage{subcaption}
\usepackage{booktabs}
\usepackage[numbers]{natbib}

\usepackage{graphics}
\usepackage{natbib}
\numberwithin{equation}{section}
\newcommand{\dm}[1]{{\displaystyle{#1}}}

\newtheorem{theorem}{\bf Theorem}[section]
\newtheorem{definition}{Definition}[section]
\newtheorem{corollary}{Corollary}[section]

\newtheorem{lemma}{Lemma}[section]
\newtheorem{remark}{Remark}[section]
\theoremstyle{remark}
\newtheorem{exam}{\bf Example}

\def \vec{\mathrm v\mathrm e \mathrm c}

\def \R{{\mathbb R}}

\def \C{{\mathbb{C}}}

\def \x{\widetilde x}
\def \y{\widetilde y}
\def \z{\widetilde z}

\def\bmatrix#1{\left[\begin{matrix}
		#1
	\end{matrix}\right]}

\def \diag{\mathrm{diag}}

\def \Df{{\mathfrak D}}
\def \A{{\mathcal A}}
\def \S{{\mathcal S }}

\def \D{{\Delta}}

\def \m{\sigma}

\def \g{{\bf g}}

\def \bdot{{\, \odot\, }}
\newcommand{\vertiii}[1]{{\left\vert\kern-0.25ex\left\vert\kern-0.25ex\left\vert #1 
    \right\vert\kern-0.25ex\right\vert\kern-0.25ex\right\vert}}
\def \rank{\mathrm{rank}}

\def\bmatrix#1{\left[ \setlength{\arrayrulewidth}{20pt}
\begin{matrix} #1 \end{matrix} \right]}
\def \noin{\noindent}
\mathchardef\mhyphen="2D
\def \R{{\mathbb R}}

\date{}

\begin{document}
	\begin{frontmatter}
    \title{Sparsity-Preserving Structured Backward Error Analysis for Double Saddle Point Problems}
    \author[inst1]{Pinki Khatun}
  \affiliation[inst1]{organization={Department of Industrial Engineering, University of Florence},%Department and Organization
			addressline={Viale Morgagni 40/44}, 
			city={Florence},
			postcode={50134}, 
			% state={Florence},
			country={Italy}}
		\ead{pinki.khatun@unifi.it, pinkikhatun.nla@gmail.com}
        \cortext[cor1]{Corresponding author.}
         \author[inst2]{Sk. Safique Ahmad\corref{cor1}}
    \affiliation[inst2]{organization={Department  of Mathematics, Indian Institute of Technology Indore},%Department and Organization
			addressline={Simrol}, 
			city={Indore},
			postcode={453552}, 
			state={Madhya Pradesh},
			country={India}}
        % \cortext[cor1]{Equal contribution.}
		\ead{safique@iiti.ac.in}
        \begin{abstract}
	Backward error (BE) analysis emerges as a powerful tool for assessing the backward stability and strong backward stability of numerical algorithms. In this paper, we explore structured BEs for a class of double saddle point problems (DSPPs) and aim to assess the strong backward stability of numerical algorithms designed to find their solutions. Our investigations preserve the inherent matrix structure and sparsity pattern in the corresponding perturbation matrices and derive explicit formulae for the structured BEs. Moreover, we provide formulae for the structure-preserving minimal perturbation matrices for which the structured BE is attained. Utilizing the relationship between the DSPP and the least squares problem with equality constraints (LSE), we derive the sparsity-preserving BE formula for LSE within our framework. Numerical experiments are performed to test the strong backward stability of various numerical algorithms.
	\end{abstract}
    
	% \noindent {\bf AMS subject classification.} 15A06, 65F10, 65F99
    \end{frontmatter}
	% \maketitle
	
%
	% \noindent {\bf Keywords.}  Backward error analysis,  Double Saddle Point Problems, Least square problems with equality constraints, Iterative algorithms, Sparse matrices.
 
	% \noindent {\bf AMS subject classification.} 15A06, 65F10, 65F99
	% \footnotetext[1]{Corresponding author.}
	% \footnotetext[2]{
	% 	Department of Mathematics, Indian Institute of Technology Indore, Khandwa Road, Indore, 452020, Madhya Pradesh, India, \texttt{safique@iiti.ac.in}.}
 %        \footnotetext[3]{Department of Industrial Engineering, University of Florence, Viale Morgagni 40/44, 50134 Florence, Italy, \texttt{pinkikhatun.nla@gmail.com}}
 %%%%%%%%%%
	%\footnotetext[3]{Research Scholar, Department of Mathematics, Indian Institute of Technology Indore, Khandwa Road, Indore, 452020, Madhya Pradesh, India, \texttt{phd2001141004@iiti.ac.in}, \texttt{pinkikhatun.nla@gmail.com}.}
  %%%%%%%%%%%%%%%%%%%%%%%%%%%%%%
% \end{proof}

% %%%%%%%%%%%%%%%%%%%%%%%%%%%%%%%%%%%%%%%
\section{Introduction}
 The double saddle point problem (DSPP) is a class of $(n+m+p)\times (n+m+p)$ block linear systems, which has attracted the attention of many researchers in recent decades for its versatile application in the community of scientific computations. For example, the DSPP arises in solving PDE-constrained optimization problems \cite{PDE-constrained2010},  liquid crystal director models \cite{LCDM2014},  Stokes-Darcy equations \cite{stokesdarcy2019}, finite element 
 discretization of Maxwell equations \cite{MaxwellEQ} and incompressible magnetohydrodynamics models \cite{MHD2019}, and so on.
We consider the general form of DSPP in the following form: %\cite{Susanne2023}:
\begin{equation}\label{SPP:EQ1}
\A {\bm w}:=\bmatrix{
A &  B_1^T & {\bf 0}\\ B_2 & -C & D_1^T\\ {\bf 0} & D_2 & E}\bmatrix{x\\y\\ z    
}=\bmatrix{f\\ g\\ h}=:{\bm d},
\end{equation}
where $A\in \R^{n\times n},$ $C\in \R^{m\times m},$ $E\in \R^{p\times p}$, $B_1,B_2\in \R^{m\times n},$ $D_1,D_2\in \R^{p\times m},$ $x,f\in \R^n,$ $y, g\in \R^{m}$ and $z,h\in \R^p.$ In most applications, the diagonal block matrices are symmetric.

Many researchers have developed efficient iterative algorithms in past years for solving the  DSPP of the form \eqref{SPP:EQ1}; for example \cite{ahmad2025robust, ahmad2026class, Susanne2023, MDS2013, Greif2023, Multiparameter2023}. To assess the stability and reliability of an approximate solution obtained using an iterative algorithm, backward error (BE) analysis is widely employed in numerical linear algebra \cite{higham2002}.  It aims to identify a closely perturbed problem (having minimum norm perturbation), ensuring that the approximate solution of the original problem coincides with the exact solution of the perturbed problem. The minimal distance between the original and perturbed problem is referred to as the BE. Furthermore, BE is used to assess the quality of an approximate solution as well as serves as a reliable and effective termination criterion when solving a problem using some iterative algorithm \cite{stopping2009}. For an approximate solution of a given problem, if the computed BE is of unit round-off error, then the corresponding algorithm is classified as backward stable \cite{higham2002}. Moreover,  the algorithm is classified as strongly backward stable \cite{strongweak, strongstab} if the perturbation matrix (with small magnitude) retains the structure of the coefficient matrix. 
% Hence, BE is also used to assess the quality of an approximate solution. %However, the considered perturbation matrix in the BE investigation does not necessarily follow the original structure of the coefficient matrix $\A$ in \eqref{SPP:EQ1}.
This leads to a natural inquiry: Whether a backward stable algorithm for solving \eqref{SPP:EQ1} exhibits strong backward stability or not?  In this paper, the notion of structured BE facilitates us in addressing the aforementioned question, where we study the BE by preserving the inherent structure of the coefficient matrix $\A$. 
 %%%%%%%%%%%%%%%%

 In various applications, such as the discretization of the Stokes equation \cite{Elman2005} and PDE-constrained optimization problems \cite{PDE-constrained2010}, the coefficient matrix of the DSPP is highly sparse. Preserving this sparsity is crucial for computational efficiency and maintaining the problem’s inherent structure. Recent studies on structured BE analysis for eigenvalue problems \cite{prince2020, PEVP2012}  and two-by-two block saddle point problem \cite{Pinki_BE, Khatun2025} have highlighted the importance of incorporating sparsity preservation into perturbation analysis. Therefore, to ensure both accuracy and efficiency in solving DSPPs, it is essential to develop optimal sparse perturbation matrices that maintain the underlying sparsity pattern.
 
Note that a DSPP of the form \eqref{SPP:EQ1} can be converted into a two-by-two block saddle point problem \cite{ZZBai2021}. For instance, if we partition the
 coefficient matrix $\A$ in \eqref{SPP:EQ1} in the following form: %However, the partition of the three-by-three block SPP in the following form
\begin{eqnarray*}
    \bmatrix{H & F_1^T\\ F_2& G}\bmatrix{p\\q}=\bmatrix{u\\v},%\quad \text{and}\quad \bmatrix{\widetilde{H} & \widetilde{F}_1^T\\ \widetilde{F}_2& \widetilde{G}}\bmatrix{\widetilde{p}\\ \widetilde{q}}=\bmatrix{\widetilde{u}\\ \widetilde{v} },
\end{eqnarray*}
where  $H=\bmatrix{A & B_1^T\\B_2 & -C},~ G=E,~ F_i=\bmatrix{{\bf 0}& D_i},$ $i=1,2,$ $p=\bmatrix{x\\y},~ q=z,~  u=\bmatrix{f\\g},$ and $ v=h.$ 
%
%$H=\bmatrix{A & B_1^T\\B_2 & -C},$ $\widetilde{H}=A,$  $F_i=\bmatrix{{\bf 0}& D_i},$ $\widetilde{F}_i=\bmatrix{B_i\\ {\bf 0}},$ $i=1,2;$ $G=E,$  $\widetilde{G}=\bmatrix{-C & D_1^T\\D_2 & E}, p=\bmatrix{x\\y}, \widetilde{p}=x, q=z, \widetilde{q}=\bmatrix{y\\z},   u=\bmatrix{f\\g},  \widetilde{u}=f,  v=h, \widetilde{v}=\bmatrix{g\\h}.$
%Then, the above partitioning leads the  DSPP \eqref{SPP:EQ1} into the standard two-by-two block $SPP$ \cite{Benzi2005}. 
Considerable research effort has been devoted to structured BEs and condition numbers for the two-by-two block saddle point problem in the past years; see \cite{PinkiGSPP, Pinki_BE, be2007wei, be2012LAA, be2017ma, BE2020BING, be2022lma}.   Further, 
 Lv \cite{threeBE2023} studied the structured BEs of the equivalent form of the DSPP \eqref{SPP:EQ1}
 given by \begin{equation}\label{SPP:EQ2}
\widehat{\A} \widehat{\bm w}:=\bmatrix{A &   {\bf 0} & B^T\\ {\bf 0} & E & D\\ -B& -D^T & C}\bmatrix{x\\z\\ y}=\bmatrix{f\\ h\\-g},
\end{equation}
with $B_1=B_2=:B,$ $D_1=D_2=:D,$ $A$ and $E$ are non-symmetric, and $C$ is symmetric. When $C={\bf 0}$ and $E={\bf 0},$ computable expressions for the structured BEs are obtained in \cite{BE2024} in three the cases: first,  $A^T=A, B_1\neq B_2$ and $D_1=D_2;$ second,  $A^T=A, B_1= B_2$ and $D_1\neq D_2;$ and third, $A^T= A, B_1\neq B_2$ and $D_1\neq D_2$. However, these studies lack the following investigations: $(a)$ the coefficient matrix $\A$ in \eqref{SPP:EQ1} is generally sparse, and the existing studies do not preserve the sparsity pattern to the perturbation matrices, $(b)$ existing research does not provide explicit formulae for the minimal perturbation matrices that preserve the inherent structures of original matrices for which an approximated solution becomes the exact solution of a nearly perturbed DSPP.

To address the aforementioned drawbacks, in this paper, we investigate structured BEs for DSPP \eqref{SPP:EQ1} by preserving sparsity pattern in three cases: $(i)$ $A^T=A,$ $B_1=B_2=:B,$ $C^T=C,$ $D_1=D_2=:D$ and $E^T=E;$   $(ii)$ $A^T= A,$ $B_1\neq B_2,$ $C^T=C,$ $D_1=D_2=:D$ and $E^T=E;$ $(iii)$ $A^T\neq A,$ $B_1=B_2=:B,$ $C^T=C,$ $D_1\neq D_2$ and $E^T=E.$

The main contributions of this study are as follows:
\begin{itemize}
    \item We investigate the structured BEs when the perturbation matrices preserve the structures mentioned in the cases $(i)$, $(ii)$ and $(iii)$, as well as preserve the sparsity patterns of the block matrices of the coefficient matrix $\A.$
    \item We derive explicit formulae for the minimal perturbation matrices for which the structured BE is attained. These perturbation matrices preserve the inherent structures of the original matrices as well as their sparsity pattern.
    \item By leveraging the connection between DSPP and least squares problems with equality constraints (LSE), we derive the sparsity-preserving BE for the LSE problem using our proposed framework.
    % \item Utilizing the derived structured BE formulae, we propose three novel termination criteria for iterative algorithms to find the effective solution of  DSPPs.
    \item Numerical experiments are performed to validate our theoretical findings and to test the backward stability and strong backward stability of numerical algorithms for solving DSPPs.
\end{itemize}
The organization of this paper is as follows. In Section \ref{sec2}, we present some notation, definitions, and preliminary results. In Sections \ref{sec3}, \ref{sec4} and \ref{sec5}, we derive explicit formulae for the structured BEs corresponding to cases $(i)$, $(ii)$ and $(iii)$, respectively.  In Section \ref{sec:LSE}, we derive the BE for the LSE problem. Section \ref{sec:numerical} includes extensive numerical experiments. Lastly, Section \ref{sec:conclusion} presents some concluding remarks.
%%%%%%%%%%%

\section{Notation, preliminaries and definitions}\label{sec2}
%%%%%%%%%%%%%%%
\subsection{Notation}
In this paper, we consistently utilize the following notation. We denote $\mathbb{R}^{m\times n}$ as the collection of all $m\times n$ real matrices and $\S_n$ as the set of all $n\times n$ symmetric matrices. The symbols $\|\cdot\|_2$ and $\|\cdot\|_F$ stand for the Euclidean and Frobenius norms, respectively.  For any matrix $X=[x_{ij}]\in \R^{m\times n},$  the symbols $X^T$ and $X^{\dagger}$ refer to the transpose and the Moore-Penrose inverse of $X,$ respectively. The notation ${\bf 1}_{m\times n}\in \R^{m\times n}$ represents the matrix with all entries are set to $1$. We set $\mu=\frac{n(n+1)}{2},$ $\sigma=\frac{m(m+1)}{2}$ and $\tau=\frac{p(p+1)}{2}.$ For  $X=[{\bf x}_1, {\bf x}_2,\ldots,{\bf x}_n]\in \R^{m\times n},$ set $\vec(X):=[{\bf{x}}^T_1,{\bf{x}}^T_2,\ldots,{\bf{x}}^T_n]^T\in \R^{mn},$ where ${\bf{x}}_i$ is the $i$-th column of  $X.$ For any symmetric matrix $X\in \S_n,$ we define its generator vector by $\vec_{\S}(X):=[{\bm x}^T_1,{\bm x}^T_2,\ldots, {\bm x}^T_n]^T\in \R^\mu,$ where ${\bm x}_1=[x_{11},x_{21},\ldots, x_{n1}]^T\in \R^n,$ ${\bm x}_2=[x_{22},x_{32},\ldots, x_{n2}]^T\in \R^{n-1}, \ldots, {\bm x}_{n-1}=[x_{(n-1)(n-1)}, x_{n(n-1)}]^T\in \R^2, {\bm x}_n=[x_{nn}]\in \R.$  The Hadamard product of $X, Y\in \R^{m\times n}$ is defined as $X\odot Y=[x_{ij}y_{ij}]\in \R^{m\times n}.$ For  $X\in \R^{m\times n}$, we define $\Theta_X:={\tt sgn}(X)=[{\tt sgn}(x_{ij})]\in \R^{m\times n},$ where
 \begin{align*}
     {\tt sgn}(x_{ij})=\left\{ \begin{array}{lcl}
     1, &\text{for}~x_{ij}\neq 0,\\
     0,  &\text{for}~ x_{ij}=0.
      \end{array}\right.
 \end{align*}
 For $x=[x_i]\in \R^n,$  $\mathfrak{D}_x$ denotes the diagonal matrix with the $i$-th diagonal entry $x_i.$
 We use the notation ${\bf 0}_{m\times n}$ to represent the zero matrix of size $m\times n$. For simplicity, we use ${\bf 0}$ when the matrix size is clear.  %$\bm{\theta}:=[\theta_1, \theta_2,\ldots, \theta_{10}]^T$, where $\theta_i$ 
Let  $\theta_i, $ $i=1,2,\ldots, 10$ be the nonnegative real numbers  with the convention that $\frac{1}{\theta_i}=0,$ whenever $\theta_i=0.$ %We define 
 % \begin{equation}
 %     \vertiii{\bmatrix{\D \A& \D{\bm d}}}_{\bm{\theta},F}:=\left\|\bmatrix{\theta_1\|\D A\|_F&  \theta_2\|\D B\|_F& \theta_3\|\D C\|_F &\theta_4\|\D D\|_F\\ \theta_5\|\D E\|_F& \theta_6\|\D f\|_2&\theta_7  \| \D g\|_2 & \theta_8\|\D h\|_2}\right\|_F.
 % \end{equation}
%%%%%%%%%%%%%%
\subsection{Preliminaries and definitions}
In this subsection, we recall the definition of unstructured BE and introduce the definitions of structured  BEs for the three cases $(i)$-$(iii).$ Furthermore, we establish two pivotal lemmas essential for deriving structured BEs. Throughout the paper, we assume that $\A$ is nonsingular.
Next, we recall the definition of unstructured BE for an approximate solution of the linear system $\A {\bm w}={\bm d}.$
\begin{definition}\label{def:UBE}
     Let $\widetilde{{\bm w}}=[\tilde{x}^{T},\tilde{y}^{T}, \tilde{z}^T]^{T}$ be an approximate solution of the DSPP \eqref{SPP:EQ1}.Then, the unstructured BE ${\bm \eta}(\tilde{{\bm w}})$ is defined as: 
 \begin{eqnarray}
 \nonumber
{\bm \eta}({\tilde{{\bm w}}})&:=&\min_{\Delta\A,\Delta {\bm d}}\left\{\left\|\bmatrix{\frac{\|\Delta \A\|_F}{\|\A\|_F}, & \frac{\|\Delta {\bm d}\|_F}{\|{\bm d}\|_F}}\right\|_2{\Big |}(\A+\Delta \A)\tilde{{\bm w}}={\bm d}+\Delta {\bm d}\right\}.
 %\end{equation}
%which has a compact formula, given by \citet{BEKKT2004} as follows: 
%\begin{equation} {\bm \eta}(\tilde{{\bm w}})
% &=&
\end{eqnarray}
\end{definition}
The unstructured BE formula is given by  Li and Liu \cite{BEKKT2004} \begin{equation}\label{BE:EQ2}
    {\bm \eta}({\tilde{{\bm w}}})=\frac{\|{\bm d}-\A \tilde{{\bm w}}\|_2}{\sqrt{\|\A\|^2_F\|\tilde{{\bm w}}\|_2^2+\|{\bm d}\|_2^2}}.
\end{equation}
 A small value of ${\bm \eta}(\widetilde{\bm w})$ indicates that the approximated solution $\widetilde{\bm w}$ is the exact solution of a slightly perturbed system $(\A +\Delta \A)\widetilde{\bm w}={\bm d}+\Delta {\bm d},$ where $\|\D \A\|_F$ and $\|\D {\bm d}\|_2$ are relatively small. That is, the corresponding iterative algorithm is backward stable.
%%%%%%%%%%%%%%%%%%%%%%%%%%%

Next, we define structured BEs for an approximate solution of the   DSPP \eqref{SPP:EQ1}.
\begin{definition}\label{def:SBE}
    Let $\widetilde{{\bm w}}=[\widetilde{x}^T, \widetilde{y}^T, \widetilde{z}^T]^T$ be an approximate solution of the  DSPP \eqref{SPP:EQ1}. Then, we define the structured BEs, denoted by ${\bm \eta}^{\mathbb{S}_i}(\x,\y,\z),$ $i=1,2,3,$  for the cases $(i)$-$(iii)$ as follows:
    \begin{align*}%\label{BE:eq21}
       & {\bm \eta}^{\mathbb{S}_1}(\widetilde{x},\widetilde{y},\widetilde{z})=\displaystyle{\min_{\left(\begin{array}{c}
       \D A,\D B,\D C,\\
     \D D, \D E, \D f,\\
      \D g, \D h
     \end{array}\right)\in\, \mathbb{S}_1}}
     \left\|\bmatrix{\theta_1\|\D A\|_F&  \theta_2\|\D B\|_F& \theta_4\|\D C\|_F\\ \theta_5\|\D D\|_F& \theta_7\|\D E\|_F& 0\\ \theta_8\|\D f\|_2&\theta_9  \| \D g\|_2 & \theta_{10}\|\D h\|_2}\right\|_F,\\
    % \vertiii{\bmatrix{\D \A& \D{\bm d}}}_{\bm{\theta},F}, \quad \text{for}\quad i=1,2,
     % 
   %  &=\left\|\bmatrix{w_1\|\D A\|_F&  w_2\|\D B\|_F& w_3\|\D C\|_F &w_4\|\D D\|_F\\ w_5\|\D E\|_F& w_6\|\D f\|_2&w_7  \| \D g\|_2 & w_8\|\D h\|_2}\right\|_F,
   %
  & {\bm \eta}^{\mathbb{S}_2}(\widetilde{x},\widetilde{y},\widetilde{z})=\displaystyle{\min_{\left(\begin{array}{c}
       \D A,\D B_1,\D B_2,\\
     \D C, \D D, \D E, \\
     \D f, \D g, \D h
     \end{array}\right)\in\, \mathbb{S}_2}}
     \left\|\bmatrix{\theta_1\|\D A\|_F&  \theta_2\|\D B_1\|_F& \theta_3\|\D B_2\|_F\\
     \theta_4\|\D C\|_F &\theta_5\|\D D\|_F& \theta_7\|\D E\|_F\\ \theta_8\|\D f\|_2&\theta_9  \| \D g\|_2 & \theta_{10}\|\D h\|_2}\right\|_F,\\
    & {\bm \eta}^{\mathbb{S}_3}(\widetilde{x},\widetilde{y},\widetilde{z})=\displaystyle{\min_{\left(\begin{array}{c}
       \D A,\D B,\D C,\\
     \D D_1, \D D_2, \D E, \\
     \D f, \D g, \D h
     \end{array}\right)\in\, \mathbb{S}_3}}
     \left\|\bmatrix{\theta_1\|\D A\|_F&  \theta_2\|\D B\|_F& \theta_4\|\D C\|_F\\
     \theta_5\|\D D_1\|_F &\theta_6\|\D D_2\|_F& \theta_7\|\D E\|_F\\ \theta_8\|\D f\|_2&\theta_9  \| \D g\|_2 & \theta_{10}\|\D h\|_2}\right\|_F,
   \end{align*}
  respectively, where 
 \begin{align}\label{s2:eq23}
     &\nonumber\mathbb{S}_1=\Bigg\{\left(\begin{array}{c}
       \D A,\D B,\D C,\\
     \D D, \D E, \D f,\\
      \D g, \D h
     \end{array}\right) {\Bigg |} \bmatrix{A+\D A & (B+\D B)^T &{\bf 0}\\ B+\D B & -(C+\D C) & (D+\D D)^T\\  {\bf 0}& D+\D D & (E+\D E)^T}\bmatrix{\widetilde{x}\\ \widetilde{y}\\ \widetilde{z}}=\bmatrix{f+\D f\\ g+\D g\\ h+\D h},\\ \nonumber
     &\hspace{4cm}\quad \D A\in \mathcal{S}_{n},\D C\in \mathcal{S}_{m}, \D E\in \mathcal{S}_{p}, \D B\in \R^{m\times n}, \D D\in \R^{p\times m}, \\ 
     &\hspace{4cm}\quad\D f\in \R^n,\D g\in \R^{m},\D h\in \R^p\Bigg\}.\\ \label{ss2:eq24}
     \end{align}
      \begin{align}
     &\nonumber\mathbb{S}_2=\Bigg\{\left(\begin{array}{c}
       \D A,\D B_1,\D B_2,\\
    \D C, \D D, \D E, \\
     \D f, \D g, \D h
     \end{array}\right) {\Bigg |} \bmatrix{A+\D A & (B_1+\D B_1)^T &{\bf 0}\\ B_2+\D B_2 & -(C+\D C) & (D+\D D)^T\\  {\bf 0}& D+\D D & (E+\D E)^T}\bmatrix{\widetilde{x}\\ \widetilde{y}\\ \widetilde{z}}=\bmatrix{f+\D f\\ g+\D g\\ h+\D h},\\ \nonumber
     &\hspace{4cm} \D A\in \S_n, \D C\in \mathcal{S}_{m}, \D E\in \mathcal{S}_{p}, \D B_1,\D B_2\in \R^{m\times n}, \D D\in \R^{p\times m}, \\
     &\hspace{4cm}\D f\in \R^n,\D g\in \R^{m},\D h\in \R^p\Bigg\}.\\ \label{s2:eq24}
     \end{align}
     \begin{align}
     &\nonumber\mathbb{S}_3=\Bigg\{\left(\begin{array}{c}
       \D A,\D B,\D C,\\
     \D D_1, \D D_2, \D E, \\
     \D f, \D g, \D h
     \end{array}\right) {\Bigg |} \bmatrix{A+\D A & (B+\D B)^T &{\bf 0}\\ B+\D B & -(C+\D C) & (D_1+\D D_1)^T\\  {\bf 0}& D_2+\D D_2 & (E+\D E)^T}\bmatrix{\widetilde{x}\\ \widetilde{y}\\ \widetilde{z}}=\bmatrix{f+\D f\\ g+\D g\\ h+\D h},\\ \nonumber
     &\hspace{3.5cm}\D A\in \R^{n\times n}, \D C\in \mathcal{S}_{m}, \D E\in \mathcal{S}_{p}, \D B\in \R^{m\times n}, \D D_1, \D D_2\in \R^{p\times m},  \\
     &\hspace{3.5cm}\D f\in \R^n,\D g\in \R^{m},\D h\in \R^p\Bigg\}.
     \end{align}
\end{definition}
%%%%%%%%%%%%%%%%%%%%
Next, we state the problem of finding structure-preserving minimal perturbation matrices for which the structured BE is attained.
 % \begin{problem}\label{Problem21}
      Further, we aim to find the minimal  perturbation matrices $\left(\begin{array}{c}
       \widehat{\D A},\widehat{\D B},\widehat{\D C},\\
     \widehat{\D D}, \widehat{\D E},\widehat{\D f },\\ \widehat{\D g}, \widehat{\D h}
     \end{array}\right)\in \mathbb{S}_1,$ $\left(\begin{array}{c}
       \widehat{\D A},\widehat{\D B_1},\widehat{\D B_2},\\
     \widehat{\D C}, \widehat{\D D}, \widehat{\D E},\\ \widehat{\D f }, \widehat{\D g}, \widehat{\D h}
     \end{array}\right)\in \mathbb{S}_2$ and   $\left(\begin{array}{c}
       \widehat{\D A},\widehat{\D B},\widehat{\D C},\\
     \widehat{\D D_1},\widehat{\D D_2 }, \widehat{\D E},\\ \widehat{\D f },\widehat{\D g}, \widehat{\D h}
     \end{array}\right)\in \mathbb{S}_3$ 
     such that
 \begin{eqnarray*}%\label{eq29}
     &{\bm \eta}^{\mathbb{S}_1}(\widetilde{x},\widetilde{y},\z) %\vertiii{\bmatrix{\widehat{\D\A},& {\widehat{\D\bm d}}}}_{\bm{\theta},F}
     =\left\|\bmatrix{\theta_1\|\widehat{\D 
 A}\|_F&  \theta_2\|\widehat{\D B}\|_F & \theta_4\|\widehat{\D C}\|_F \\ \theta_5\|\widehat{\D D}\|_F &\theta_7\|\widehat{\D E}\|_F& 0\\ \theta_8\|\widehat{\D f}\|_2& \theta_9  \| \widehat{\D g}\|_2 & \theta_{10}\|\widehat{\D h}\|_F}\right\|_F, \\
 %%%
 &{\bm \eta}^{\mathbb{S}_2}(\widetilde{x},\widetilde{y},\z)=\left\|\bmatrix{\theta_1\|\widehat{\D 
 A}\|_F&  \theta_2\|\widehat{\D B_1}\|_F & \theta_3\|\widehat{\D B_2}\|_F \\ \theta_4\|\widehat{\D C}\|_F &\theta_5\|\widehat{\D D}\|_F& \theta_7\|\widehat{\D E}\|_F\\ \theta_8\|\widehat{\D f}\|_2& \theta_9  \| \widehat{\D g}\|_2 & \theta_{10}\|\widehat{\D h}\|_F}\right\|_F,%\\
 %%%
 \end{eqnarray*}
 \begin{eqnarray*}
 & {\bm \eta}^{\mathbb{S}_3}(\widetilde{x},\widetilde{y},\z) =\left\|\bmatrix{\theta_1\|\widehat{\D 
 A}\|_F&  \theta_2\|\widehat{\D B}\|_F & \theta_4\|\widehat{\D C}\|_F \\ \theta_5\|\widehat{\D D_1}\|_F &\theta_6\|\widehat{\D D_2}\|_F&\theta_7\|\widehat{\D E}\|_F\\ \theta_8\|\widehat{\D f}\|_2& \theta_9  \| \widehat{\D g}\|_2 & \theta_{10}\|\widehat{\D h}\|_F}\right\|_F,
 \end{eqnarray*}
 respectively.
 % \end{problem}
%%%%%
\begin{remark}
    When $\theta_i = 0$ for any given $i$ ($i=1,2,\ldots,10$), it indicates that the corresponding block matrix has no perturbation.
\end{remark}
\begin{remark}\label{REM1}
% Many researchers have studied the structured BE by preserving the sparsity pattern; see \cite{ Khatun2025,prince2021, Pinki_BE, zhang2012structured}.
Our specific interest lies in investigating structured BEs while the perturbation matrices preserve the sparsity pattern of the original matrices. To perform this investigation, we substitute the perturbation matrices $\Delta A$, $\Delta B_1$, $\D B_2,$ $\Delta C$, $\Delta D_1$ $\D D_2$ and $\Delta E$ by $\Delta A \odot \Theta_A,$ $\Delta B_1 \odot \Theta_{B_1},$ $\Delta B_2 \odot \Theta_{B_2},$ $\Delta C \odot \Theta_C,$ $\Delta D_1 \odot \Theta_{D_1},$ $\Delta D_2 \odot \Theta_{D_2}$ and $\Delta E \odot \Theta_E,$ respectively. Within this framework, we denote the structured BEs as ${\bm \eta}_{\textbf{sps}}^{\mathbb{S}_i}(\x,\y,\z),$ $i=1,2,3.$ Moreover, the minimal perturbation matrices are denoted by $ \widehat{\D A}_{\tt sps},$ $\widehat{\D B_1}_{\tt sps},$ $ \widehat{\D B_2}_{\tt sps},$ $\widehat{\D C}_{\tt sps},$ $ \widehat{\D D_1}_{\tt sps},$ $ \widehat{\D D_2}_{\tt sps},$  $ \widehat{\D E}_{\tt sps},$ $\widehat{\D f}_{\tt sps},$ $ \widehat{\D g}_{\tt sps}$ and  $\widehat{\D h}_{\tt sps}.$
\end{remark}
 When the structured BEs ${\bm \eta}^{\mathbb{S}_i}(\widetilde{x},\widetilde{y},\z)$  and ${\bm \eta}_{\tt sps}^{\mathbb{S}_i}(\widetilde{x},\widetilde{y},\z)$  are around an order of unit round-off error, then the approximate solution $\widetilde{\bm w}=[\widetilde{x}^T,\widetilde{y}^T,\z^T]^T$ becomes an exact solution of nearly perturbed structure-preserving DSPP of the form \eqref{SPP:EQ1}. Thus, the corresponding algorithm is referred to as strongly backward stable.
 To obtain the structured BEs formulae, the following lemmas play a pivotal role.
\begin{lemma}\label{Lemma1}
     Let $A, H\in\mathcal{S}_n$ with  generator vectors $\vec_{\mathcal{S}}(A)=[\bm{a}^T_1,\bm{a}^T_2\ldots,\bm{a}^T_n]^T$ and \(\vec_{\mathcal{S}}(H)=[\bm{h}^T_1,\bm{h}^T_2\ldots,\bm{h}^T_n]^T,\) respectively. Suppose $x=[x_1,\ldots,x_n]^T\in \R^n,$   $y=[y_1,\ldots,y_n]^T\in \R^n$ and ${\bm d}=[d_1,\ldots,d_n]^T\in \R^n.$  Then $(A\,\odot\,\Theta_{H})x=\bm{d}$  can be expressed as:
     \begin{align*}
       \mathcal{K}_{x}\Phi_H  \vec_{\mathcal{S}}(A\,\odot\,\Theta_{H})={\bm d},
     \end{align*}
where  $\Phi_H=\diag(\vec_{\mathcal{S}}(\Theta_H)),$  $\mathcal{K}_x=\bmatrix{K^1_x& K^2_x &\cdots &K^n_x}\in \R^{n\times \mu}$ and $K^i_x\in \R^{n\times (n-i+1)}$ are given by
\begin{align*}
  %%%%
  &K^1_x=\bmatrix{x_1& x_2&\cdots & \cdots& x_n\\
    0 &x_1 &0& \cdots &0\\
    0&0&x_1&\cdots&0\\
    \vdots& \vdots &\ddots &\ddots& \vdots\\
    0&\cdots &\cdots& 0& x_1
    },~ K^2_x=\bmatrix{0 & 0&\cdots&\cdots &0\\
    x_2&x_3& \cdots & \cdots& x_n\\
    0 &x_2 &0& \cdots &0\\
    0&0&x_2&\cdots&0\\
    \vdots& \vdots &\ddots &\ddots& \vdots\\
    0&\cdots &\cdots& 0& x_2
    },
    % \end{align*}
    % \begin{align*}
   \ldots,~ K^n_x=\bmatrix{0\\0\\ \vdots\\ 0\\x_n}.
\end{align*}
\end{lemma}
{\em Proof.}  Since $H\in \S_n,$ we have $\Theta_H\in \S_n.$ Let $\vec_{\S}(A)=[{\bm a}_1,{\bm a}_2,\ldots {\bm a}_n]^T$ and $\vec_{\mathcal{S}}(\Theta_H)=[\Theta_{{\bm h}_1}^T, \Theta_{{\bm h}_2}^T,\ldots, \Theta_{{\bm h}_n}^T ]^T,$ where $\Theta_{{\bm h}_i}=[{\tt sgn}(h_{ii}), {\tt sgn}(h_{i(i+1)}),\ldots, {\tt sgn}(h_{in})]^T\in \R^{n-i+1}.$  
 Then, $(A\,\odot\,\Theta_{H})x$ can be equivalently expressed as 
\begin{eqnarray}
 \nonumber (A\,\odot\,\Theta_{H}) x &=&\bmatrix{a_{11}\, {\tt sgn}(h_{11})x_1+ a_{21}\, {\tt sgn}(h_{21})x_2+\cdots+ a_{n1}\, {\tt sgn}(h_{n1})x_n\\
       a_{21}\, {\tt sgn}(h_{21})x_1+ a_{22}\, {\tt sgn}(h_{22})x_2+\cdots+ a_{n2}\, {\tt sgn}(h_{n2})x_n\\
        a_{31}\, {\tt sgn}(h_{31})x_1+ a_{32}\, {\tt sgn}(h_{32})x_2+\cdots+ a_{n3}\, {\tt sgn}(h_{n3})x_n\\
       \vdots \\
       % a_{(n-1)1}\, {\tt sgn}(h_{(n-1)1})x_1+ a_{(n-1)2}\, {\tt sgn}(h_{(n-1)2})x_2+\cdots+ a_{nn}\, {\tt sgn}(h_{nn})x_n\\
      a_{n1}\, {\tt sgn}(h_{n1})x_1+ a_{n2}\, {\tt sgn}(h_{n2})x_2+\cdots+ a_{nn}\, {\tt sgn}(h_{nn})x_n}\\[2ex] \nonumber
      %%%%%%
      &=&\bmatrix{a_{11}\, {\tt sgn}(h_{11})x_1+ a_{21}\, {\tt sgn}(h_{21})x_2+\cdots+ a_{1n}\, {\tt sgn}(h_{n1})x_n\\
       a_{21}{\tt sgn}(h_{21})x_1\\
       \vdots \\
      a_{n1}\, {\tt sgn}(h_{n1})x_1}\\\nonumber
      &+&\bmatrix{0\\
        a_{22}\, {\tt sgn}(h_{22})x_2+\cdots+ a_{n2}\, {\tt sgn}(h_{n2})x_n\\
        a_{32}\, {\tt sgn}(h_{32})x_2\\
       \vdots \\
      a_{n2}\, {\tt sgn}(h_{n2})x_2}\\  \nonumber
      \end{eqnarray}
      \begin{eqnarray}  \nonumber
      %%%%%%%%%%%%%%
      &+&\bmatrix{0\\0\\
        a_{33}\, {\tt sgn}(h_{33})x_3+\cdots+ a_{n3} {\tt sgn}(h_{n3})x_n\\
       \vdots \\
      a_{n3}\, {\tt sgn}(h_{n3})x_3}+\cdots+\bmatrix{0\\
       0\\
        \vdots \\
        % a_{(n-1)(n-1)}\, {\tt sgn}(h_{(n-1)(n-1)})x_{n-1}+ a_{(n-1)n}\, {\tt sgn}(h_{32})x_{n}\\
      % a_{n2}\, {\tt sgn}(h_{n2})x_{n-1}
       0\\
       a_{nn}\, {\tt sgn}(h_{nn})x_n}\\  \nonumber
       %%%%%%%%%%%%
       &=&K^1_x\mathfrak{D}_{ \Theta_{{\bm h}_1}}({\bm a}_1\odot \Theta_{{\bm h}_1})+K^2_x \mathfrak{D}_{ \Theta_{{\bm h}_2}} ({\bm a}_2\odot \Theta_{{\bm h}_2}) + \cdots+ K^n_x \mathfrak{D}_{ \Theta_{ {\bm h}_n}} ({\bm a}_n\odot \Theta_{{\bm h}_n})\\ \nonumber     
       % \end{eqnarray}
       %  \begin{eqnarray}
       &= &\bmatrix{K^1_x\mathfrak{D}_{ \Theta_{{\bm h}_1}} &K^2_x\mathfrak{D}_{\Theta_{{\bm h}_2}} & \cdots & K^n_x \mathfrak{D}_{\Theta_{{\bm h}_n}}} \bmatrix{{\bm a}_1\odot \Theta_{{\bm h}_1}\\{\bm a}_2\odot \Theta_{{\bm h}_n} \\ \vdots \\{\bm a}_n\odot \Theta_{{\bm h}_n}}\\
       &=&\mathcal{K}_{x} \Phi_H \vec_{\mathcal{S}}(A\,\odot\,\Theta_{H}).
\end{eqnarray}
Therefore, the proof is completed. $\square$ 
%%%%%%%%%%%%%%%%%%%%
% \begin{remark}
%      When the matrix $\Theta_H={\bf 1}_{n\times n},$ from Lemma \ref{Lemma1}, the linear system $Ax={\bm b},$  where $A\in \mathcal{S}_n,$ can be written  as ${\mathcal{K} }_x\vec_{\S}(A)={\bm b}$.
%  \end{remark}
 %%%%%%%%%%%%%%%%%%%%%%%%
\begin{lemma}\label{Lemma2}
       Let $A,B,H\in \R^{m\times n}$ be  three given matrices. Suppose that $x=[x_1,\ldots,x_n]^T\in \R^n$, $y=[y_1,\ldots,y_m]^T\in \R^m,$  ${\bm d}^1\in \R^m$ and ${\bm d}^2\in \R^n.$ Then $(A\,\odot\,\Theta_{H})x={\bm d}^1$ and $(B\bdot \Theta_H)^Ty={\bm d}^2$  are equivalent to
     $${\mathcal{M}}^{m}_x\mathfrak{D}_{\vec(\Theta_H)}\vec(A\,\odot\,\Theta_{H})={\bm d}^1 \quad \mbox{and } \quad { \mathcal{N}}^{n}_y\mathfrak{D}_{\vec(\Theta_H)}\vec(B\,\odot\,\Theta_{H})={\bm d}^2,$$ 
  respectively, where ${\mathcal{M}}^{m}_x=\bmatrix{x_1I_m&x_2I_m&\cdots&x_nI_m}\in \R^{m\times {mn}}$ and
  \begin{align}
{\mathcal{N}}^{n}_{y}=\bmatrix{y^T& {\bf 0}&\cdots&\cdots&{\bf 0}\\
         {\bf 0} & y^T&  {\bf 0}&\cdots&{\bf 0} \\
        \vdots& \ddots&\ddots&\ddots&\vdots\\
         \vdots& &\ddots & \ddots&{\bf 0}\\
        {\bf 0}&\cdots &\cdots &{\bf 0}&y^T}\in \R^{n\times {mn}}.
     \end{align}
    % Here, ${\bf 0}$ denotes the zero vector of order $1\times n.$
 \end{lemma}
 %%%%%%%%%%%%%%%%
 {\em Proof.} The proof proceeds using an analogous method to the proof of Lemma \ref{Lemma1}. $\square$
 %%%%%%%%%%%%%
 \begin{lemma}\label{sec2:lemma}\cite{golubvan}
    Let $A \in \mathbb{C}^{m \times n}$ and $b \in \mathbb{C}^m$. The system of linear equations $Ax = b$ is consistent if and only if $AA^\dagger  b = b$. Additionally, the minimum norm least-squares solution of this system is represented by $A^\dagger b$.
\end{lemma}

\section{Derivation of structured BEs for case $(i)$}\label{sec3}
%%%%%%%%%%%%%%
In this section, we discuss the structured BEs for the  DSPP \eqref{SPP:EQ1} for the case $(i),$ i.e., $A\in \S_n, C\in \S_m,$  $E\in \S_p,$ $B_1=B_2=:B$ and $D_1=D_2=:D,$ and perturbation matrices belongs to set $\mathbb{S}_1.$ Prior to that, we construct   the diagonal matrix $\mathfrak{D}_{\S_{n}}\in \R^{\mu\times \mu},$ where
 $$\left\{ \begin{array}{lcl}
		\mathfrak{D}_{\S_{n}}(k,k)=1, & \mbox{for} & k=\frac{(2n-(i-2))(i-1)}{2}+1,\, i=1,2,\ldots,n, \\ \mathfrak{D}_{\S_{n}}(k,k)=\sqrt{2},&  &\text{otherwise}.
	\end{array}\right.$$
The matrix $\mathfrak{D}_{\S_{n}}$ has the property, $\|A\|_F=\|\mathfrak{D}_{\S_{n}}\vec_{\mathcal{S}}(A)\|_2.$ Further, we introduce the following notation:
 \begin{align}
     &\Phi_A=\diag(\vec_{\S}(\Theta_A)),~ \Phi_B=\diag(\vec(\Theta_B)),~ \Phi_C=\diag(\vec_{\S}(\Theta_C)),\\
     &\Phi_D=\diag(\vec(\Theta_D)),~\Phi_E=\diag(\vec_{\S}(\Theta_E)),
 \end{align}
 and $$\mathcal{I}=\bmatrix{ -\frac{1}{\theta_8}I_n& {\bf 0} & {\bf 0}\\  {\bf 0} & -\frac{1}{\theta_9} I_m & {\bf 0} \\{\bf 0} & {\bf 0}&  -\frac{1}{\theta_{10}}I_p}.$$

%%%%%%%%%%%%%%%%%%%%
\begin{theorem}\label{theorem1}
     Let $[\x^T,\, \y^T,\z^T]^T$ be  an approximate solution of the  DSPP \eqref{SPP:EQ1} with $A\in \S_n,$ $C\in \S_m,$   $E\in \S_p,$ and $\theta_8, \theta_9, \theta_{10}\neq 0.$  Then, we have 
\begin{align}
    {\bm \eta}^{\mathbb{S}_1}_{\bf sps}(\x,\y,\z)=\left\|\mathcal{J}_{\mathbb{S}_1}^{T}(\mathcal{J}_{\mathbb{S}_1}\mathcal{J}_{\mathbb{S}_1}^T)^{-1}R_{\bm d}\right\|_2,
\end{align}
where $\mathcal{J}_{\mathbb{S}_1}=[\widetilde{\mathcal{J}}_{\mathbb{S}_1}~~\mathcal{I}]\in \R^{(n+m+p)\times \bm{l}}$ and  $ \widetilde{\mathcal{J}}_{\mathbb{S}_1}$ is given by
\begin{align*}
 \widetilde{\mathcal{J}}_{\mathbb{S}_1}=\bmatrix{ \frac{1}{\theta_1}{\mathcal{K}}_{\x} \Phi_A \mathfrak{D}_{\S_{n}}^{-1}& \frac{1}{\theta_2}{\mathcal{N}}^{n}_{\y}\Phi_B& {\bf 0} & {\bf 0} & {\bf 0} \\
 {\bf 0}& \frac{1}{\theta_2}M^m_{\x}\Phi_B & -\frac{1}{\theta_4}{\mathcal{K}}_{\y}\Phi_C\mathfrak{D}_{\S_{m}}^{-1} &\frac{1}{\theta_5}N^m_{\z}\Phi_D & {\bf 0} \\
 {\bf 0} & {\bf 0} &{\bf 0}& \frac{1}{\theta_5}{\mathcal{M}}^{p}_{\y}\Phi_D & \frac{1}{\theta_7}{\mathcal{K}}_{\z} \Phi_E\mathfrak{D}_{\S_{p}}^{-1}},
\end{align*}
$R_{\bm d}=[R_f^T,\, R_g^T,\, R_h^T]^T,$ $R_f=f-A\x-B^T\y,$ $R_g=g-B\x+C\y-D^T\z,$  $R_h=h-D\y-E\z,$ and $\bm{l}={\mu}+{\sigma}+{\tau}+mn+mp+m+n+p.$

\noin The minimal perturbation matrices  are given by the following generating vectors:
 \begin{eqnarray*}
  \vec_{\S}(  \widehat{\D A}_{\tt sps})&=&\theta_1^{-1}\Df_{\S_n}^{-1}\bmatrix{I_\mu & \bf 0}\mathcal{J}_{\mathbb{S}_1}^{T}(\mathcal{J}_{\mathbb{S}_1}\mathcal{J}_{\mathbb{S}_1}^T)^{-1}R_{\bm d},\\
   \vec( \widehat{\D B}_{\tt sps})&=&\theta_2^{-1}\bmatrix{\bf 0& I_{mn} & \bf 0}\mathcal{J}_{\mathbb{S}_1}^{T}(\mathcal{J}_{\mathbb{S}_1}\mathcal{J}_{\mathbb{S}_1}^T)^{-1}R_{\bm d},\\
     \vec_{\S}(\widehat{\D C}_{\tt sps})&=&\theta_4^{-1}\Df_{\S_m}^{-1}\bmatrix{\bf 0 & I_\sigma & \bf 0}\mathcal{J}_{\mathbb{S}_1}^{T}(\mathcal{J}_{\mathbb{S}_1}\mathcal{J}_{\mathbb{S}_1}^T)^{-1}R_{\bm d},\\
    \vec( \widehat{\D D}_{\tt sps})&=&\theta_5^{-1}\bmatrix{\bf 0& I_{mp} & \bf 0}\mathcal{J}_{\mathbb{S}_1}^{T}(\mathcal{J}_{\mathbb{S}_1}\mathcal{J}_{\mathbb{S}_1}^T)^{-1}R_{\bm d},\\
    \vec_{\S}( \widehat{\D E}_{\tt sps})&=&\theta_7^{-1}\Df_{\S_p}^{-1}\bmatrix{\bf 0&I_\tau & \bf 0}\mathcal{J}_{\mathbb{S}_1}^{T}(\mathcal{J}_{\mathbb{S}_1}\mathcal{J}_{\mathbb{S}_1}^T)^{-1}R_{\bm d},\\
       \widehat{\D f}_{\tt sps}&=&\theta_8^{-1}\bmatrix{\bf 0&I_{ n} & \bf 0}\mathcal{J}_{\mathbb{S}_1}^{T}(\mathcal{J}_{\mathbb{S}_1}\mathcal{J}_{\mathbb{S}_1}^T)^{-1}R_{\bm d},\\
         \widehat{\D g}_{\tt sps}&=&\theta_9^{-1}\bmatrix{\bf 0&I_{m} & \bf 0}\mathcal{J}_{\mathbb{S}_1}^{T}(\mathcal{J}_{\mathbb{S}_1}\mathcal{J}_{\mathbb{S}_1}^T)^{-1}R_{\bm d},\\
           \widehat{\D h}_{\tt sps}&=&\theta_{10}^{-1}\bmatrix{\bf 0 & I_{p}}\mathcal{J}_{\mathbb{S}_1}^{T}(\mathcal{J}_{\mathbb{S}_1}\mathcal{J}_{\mathbb{S}_1}^T)^{-1}R_{\bm d}.
 \end{eqnarray*}
\end{theorem}
%%%%%%%%%%%
\proof For the approximate solution $[\widetilde{x}^T,\widetilde{y}^T,\z^T]^T,$ we need to construct  perturbations $
       \D A\in \S_n,$ $\D B\in \R^{m\times n},$ $
     \D C\in \S_m,$ $\D D\in \R^{p\times m},$   $\D E\in \S_p,$  which maintain the sparsity pattern of $A, B, C, D, E, $ respectively, and the perturbations $ \D f\in \R^n,$ $ \D g\in \R^m,$ and $ \D h\in \R^p$. By \eqref{s2:eq23}, $\left(\begin{array}{c}
       \D A,\D B,\D C,\\
     \D D, \D E, \D f,\\
      \D g, \D h
     \end{array}\right)\in\, \mathbb{S}_1$ if and only if $\D A, \D B, \D C, \D D, \D E, \D f, \D g$ and $\D h$ satisfy
 \begin{eqnarray}\label{theorem1:eq33}
      \left.   \begin{array}{lcl}
       \D A\x+ \D B^{T}\y -\D f= R_f, &  \\
       \D B\x-\D C\y+\D D^{T}\z- \D g= R_g,& \\
       \D D\y+\D E\z-\D h=R_h,
    \end{array}\right\}
    \end{eqnarray}
    and $\D A\in \S_n,$ $\D C\in \S_m,$ $\D E\in \S_p.$
 To maintain the sparsity pattern of $A,B, C, D$ and $E$ to the perturbation matrices, we replace  $\D A,\D B, \D C,\D D$ and
     $\D E $  by $\D A\bdot \Theta_A,\D B\bdot \Theta_B,$ $ \D C\bdot \Theta_C, \D D\bdot \Theta_D$ and $\D E\bdot \Theta_E,$ respectively. Thus \eqref{theorem1:eq33} can be equivalently reformulated as:
     \begin{align}\label{theorem1:eq34}
         &\theta_1^{-1}\theta_1(\D A\bdot \Theta_A)\widetilde{x}+\theta_2^{-1}\theta_2(\D B\bdot \Theta_B)^T \y-\theta_8^{-1}\theta_8\D f=R_f,\\ \label{theorem1:eq35}
         &\theta_2^{-1}\theta_2(\D B\bdot \Theta_B)\widetilde{x}-\theta_4^{-1}\theta_4(\D C\bdot \Theta_C) \y+\theta_5^{-1}\theta_5(\D D\bdot \Theta_D)^T\z-\theta_9^{-1}\theta_9\D g=R_g,\\ \label{theorem1:eq36}
        &\theta_5^{-1}\theta_5(\D D\bdot \Theta_D)\widetilde{y}+\theta_7^{-1}\theta_7(\D E\bdot \Theta_E)^T \z-\theta_{10}^{-1}\theta_{10}\D h=R_h.
     \end{align}
Using Lemma \ref{Lemma1} in  \eqref{theorem1:eq34}, we get
\begin{align}\label{theorem1:eq37}
    \theta_1^{-1}{\mathcal{K}}_{\x} \Phi_A\theta_1\vec_{\S}(\D A\odot \Theta_{A})+  \theta_2^{-1}{\mathcal{N}}^n_{\y} \Phi_B\theta_2\vec(\D B\odot \Theta_B)- \theta_8^{-1}\theta_8\D f=R_f.
\end{align}
Further, (\ref{theorem1:eq37})  can be express as 
\begin{align}\label{theorem1:eq38}
    \theta_1^{-1}{\mathcal{K}}_{\x} \Phi_A \mathfrak{D}_{\S_{n}}^{-1} \mathfrak{D}_{\S_{n}}\theta_1\vec_{\S}(\D A\odot \Theta_{A})+  \theta_2^{-1}{\mathcal{N}}^n_{\y} \Phi_B\theta_2\vec(\D B\odot \Theta_B)- \theta_8^{-1}\theta_8\D f=R_f.
\end{align}
  Equivalently, (\ref{theorem1:eq38}) can be written as follows:
\begin{equation}\label{eq38}
    \mathcal{J}^1_{\mathbb{S}_1}\D X=R_f,
\end{equation}
where $$\mathcal{J}^1_{\mathbb{S}_1}=\bmatrix{ \theta_1^{-1}{\mathcal{K}}_{\x} \Phi_A \mathfrak{D}_{\S_{n}}^{-1}& \theta_2^{-1}{\mathcal{N}}^n_{\y}\Phi_B& {\bf 0} & {\bf 0} & {\bf 0} & -\theta^{-1}_8I_n& {\bf 0} & {\bf 0}}\in \R^{n\times {\bm l}}$$ and
\begin{equation}\label{theorem:eq317}
    \D X=\bmatrix{\theta_1\mathfrak{D}_{\S_{n}}\vec_{\S}(\D A\odot \Theta_{A})\\ \theta_2\vec(\D B\odot \Theta_B)\\ \theta_4\mathfrak{D}_{\S_{m}}\vec_{\S}(\D C\odot \Theta_{C})\\ \theta_5\vec(\D D\odot \Theta_D) \\ \theta_7\mathfrak{D}_{\S_{p}}\vec_{\S}(\D E\odot \Theta_{E})\\ \theta_8\D f\\ \theta_9\D g\\ \theta_{10}\D h}\in \R^{\bm l} .
\end{equation}

\noin Similarly, using Lemma \ref{Lemma1} to \eqref{theorem1:eq35} and \eqref{theorem1:eq36}, we obtain
\begin{equation}\label{eq39}
    \mathcal{J}^2_{\mathbb{S}_1}\D X= R_g \quad \text{and} \quad  \mathcal{J}^3_{\mathbb{S}_1}\D X= R_h,
\end{equation}
 where $\mathcal{J}^2_{\mathbb{S}_1}\in \R^{m\times {\bm l}}$ and $\mathcal{J}^3_{\mathbb{S}_1}\in \R^{p\times {\bm l}}$ are given by 
 \begin{align}
     \mathcal{J}^2_{\mathbb{S}_1}=\bmatrix{{\bf 0}& \theta_2^{-1}{\mathcal{M}}^m_{\x}\Phi_B & -\theta_4^{-1}{\mathcal{K}}_{\y}  \Phi_C \Df_{\S_m}^{-1}&\theta_5^{-1}{\mathcal{N}}^m_{\z} \Phi_D& {\bf 0} &  {\bf 0} & -\theta^{-1}_9 I_m & {\bf 0}}
 \end{align}
 and 
 \begin{equation}
     \mathcal{J}^3_{\mathbb{S}_1}=\bmatrix{{\bf 0} & {\bf 0} &{\bf 0}& \theta_5^{-1}{\mathcal{M}}^p_{\y}\Phi_D & \theta_7^{-1}{\mathcal{K}}_{\z} \Df_{\S_p}^{-1} \Phi_E & {\bf 0} & {\bf 0}&  -\theta^{-1}_{10}I_p},
 \end{equation}
respectively.
Combining \eqref{eq38} and \eqref{eq39}, we obtain the following equivalent system 
\begin{equation}\label{eq312}
    \mathcal{J}_{\mathbb{S}_1}\D X=R_{\bm d}.
\end{equation}
\noin Clearly, for $\theta_8,\theta_9,\theta_{10}\neq 0,$ $\mathcal{J}_{\mathbb{S}_1}$ has full row rank. Therefore, by Lemma \ref{sec2:lemma}, the minimum norm solution of \eqref{eq312} is given by
\begin{equation}\label{eq313}
    \D X_{\min}=\mathcal{J}_{\mathbb{S}_1}^{\dagger}\bmatrix{R_f \\R_g \\R_h} = \mathcal{J}_{\mathbb{S}_1}^T(\mathcal{J}_{\mathbb{S}_1}\mathcal{J}_{\mathbb{S}_1}^T)^{-1} R_{\bm d}.
    \end{equation}
On the other hand, the minimization problem in Definition \ref{def:SBE} can be reformulated as:
\begin{align}\label{eq314}
  \nonumber [{\bm \eta}^{\mathbb{S}_1}_{\tt sps}(\widetilde{x},\widetilde{y},\z)]^2&=\displaystyle{\min}\Bigg\{\theta_1^2\|\D A\bdot \Theta_A\|^2_F+\theta_2^2\|\D B\bdot \Theta_B\|_F^2 +\theta_4^2\|\D C\bdot \Theta_C\|_F^2 +\theta_5^2\|\D D\bdot \Theta_D\|_F^2\\ \nonumber
        &\quad \quad \quad \quad + \theta_7^2\|\D E\bdot \Theta_E\|_F^2 +\theta_8^2\|\D f\|_2^2+\theta_9^2\| \D g\|_2^2+ \theta_{10}^2\| \D h\|_2^2 \,\, {\Big |}\\ \nonumber
        &\quad \quad \quad \quad{\left(\begin{array}{c}
       \D A\bdot \Theta_A,\D B\bdot \Theta_B, \D C\bdot \Theta_C,\\ \nonumber
     \D D\bdot \Theta_D, \D E\bdot \Theta_E, \D f,\\
     \D g, \D h
     \end{array}\right)\in \mathbb{S}_1}\Bigg\}\\ \nonumber
     %%%%%
     &   =\min\Bigg\{\theta_1^2\|\mathfrak{D}_{\S_n}\vec_{\S}(\D A\odot \Theta_A)\|_2^2 +\theta_2^2\|\vec(\D B\odot \Theta_B)\|_2^2 \\ \nonumber
     %%%%%
        &\quad \quad \quad + \theta_4^2\|\mathfrak{D}_{\S_m}\vec_{\S}(\D C\odot \Theta_C)\|_2^2 +\theta_5^2\|\vec(\D D\odot \Theta_D)\|_2^2\\ \nonumber
        %%%
        &\quad \quad \quad  + \theta_7^2\|\mathfrak{D}_{\S_p}\vec_{\S}(\D E\bdot \Theta_E)\|_2^2 +\theta_8^2\|\D f\|_2^2+\theta_9^2\| \D g\|_2^2 +\theta_{10}^2\| \D h\|_2^2 \,\, {\Big |}\,\\
        &\quad \quad \quad\quad\mathcal{J}_{\mathbb{S}_1}\D X=R_{\bm d}\Bigg\}          \\
        &=\min \left\{\|\D X\|_2^2\,{\Big |}\,\mathcal{J}_{\mathbb{S}_1}\D X=R_{\bm d}\right\}=\|\D X_{\min}\|^2_2.
\end{align}
Consequently, substituting \eqref{eq313} into \eqref{eq314}, we obtain 
\begin{equation*}
{\bm \eta}^{\mathbb{S}_1}_{\tt sps}(\x,\y,\z) =\left\|\mathcal{J}_{\mathbb{S}_1}^T(\mathcal{J}_{\mathbb{S}_1}\mathcal{J}_{\mathbb{S}_1}^T)^{-1} R_{\bm d}\right\|_2.
\end{equation*}
From \eqref{theorem:eq317}, we have $\theta_1\Df_{\S_n}\vec_{\S}(\D A\odot \Theta_A)=\bmatrix{I_\mu & {\bf 0}}\D X.$ Therefore, the generating vector for the minimal perturbation matrix $\widehat{\D A}_{\tt sps}$ which also preserves the sparsity pattern  is given by $$\vec_{\S}(\widehat{\D A}_{\tt sps})=\theta_1^{-1}\Df_{\S_n}^{-1}\bmatrix{I_\mu &{ \bf 0}}\D X_{\min}.$$
Similarly, the generating vectors for other minimal perturbation matrices can be obtained. Hence, the proof is completed. $\square$
\begin{corollary}\label{coro1:theorem31}
   Suppose the  approximate solution of the  DSPP \eqref{SPP:EQ1} with $A\in \S_n,$ $C\in \S_m,$   $E\in \S_p,$ and $\theta_8, \theta_9, \theta_{10}\neq 0$ is $[\x^T,\, \y^T,\z^T]^T.$  Then, we have 
\begin{align}
    {\bm \eta}^{\mathbb{S}_1}(\x,\y,\z)=\left\|\widehat{\mathcal{J}}_{\mathbb{S}_1}^{T}(\widehat{\mathcal{J}}_{\mathbb{S}_1}\widehat{\mathcal{J}}_{\mathbb{S}_1}^T)^{-1}R_{\bm d}\right\|_2,
\end{align}
where $ \widehat{\mathcal{J}}_{\mathbb{S}_1}\in \R^{(n+m+p)\times \bm{l}}$ is given by
{\footnotesize
\begin{align*}
\widehat{ \mathcal{J}}_{\mathbb{S}_1}=\bmatrix{ \frac{1}{\theta_1}{\mathcal{K}}_{\x} \mathfrak{D}_{\S_{n}}^{-1}& \frac{1}{\theta_2}{\mathcal{N}}^{n}_{\y}& {\bf 0} & {\bf 0} & {\bf 0} & -\frac{1}{\theta_8}I_n& {\bf 0} & {\bf 0}\\
 {\bf 0}& \frac{1}{\theta_2}{\mathcal{M}}^m_{\x} & -\frac{1}{\theta_4}{\mathcal{K}}_{\y}\mathfrak{D}_{\S_{m}}^{-1} &\frac{1}{\theta_5}{\mathcal{N}}^m_{\z} & {\bf 0} &  {\bf 0} & -\frac{1}{\theta_9} I_m & {\bf 0} \\
 {\bf 0} & {\bf 0} &{\bf 0}& \frac{1}{\theta_5}{\mathcal{M}}^{p}_{\y} & \frac{1}{\theta_7}{\mathcal{K}}_{\z} \mathfrak{D}_{\S_{p}}^{-1} & {\bf 0} & {\bf 0}&  -\frac{1}{\theta_{10}}I_p }.
\end{align*}}
% $R_f=f-A\x-B^T\y,$ $R_g=g-B\x+C\y-D^T\z,$  $R_h=h-D\y-E\z,$ and $\bm{l}={\mu}+{\sigma}+{\tau}+mn+mp+m+n+p.$  
\end{corollary}
\proof Since we are not preserving the sparsity pattern, the proof follows by considering $\Theta_A={\bf 1}_{n\times n},$ $\Theta_B={\bf 1}_{m\times n},$ $\Theta_C={\bf 1}_{m\times m},$ $\Theta_D={\bf 1}_{p\times m},$ and  $\Theta_E={\bf 1}_{p\times p}$ in Theorem \ref{theorem1}.  $\square$

\begin{remark}
The structure-preserving minimal perturbation matrices $\widehat{\D A},$ $\widehat{\D B},$ $\widehat{\D C},$ $\widehat{\D D},$ $\widehat{\D E},$ $\widehat{\D f},$ $\widehat{\D g},$ and $\widehat{\D h},$  for which ${\bm \eta}^{\mathbb{S}_1}(\x,\y,\z)$ is attained are given by formulae presented in Theorem \ref{theorem1} with $\mathcal{J}_{\mathbb{S}_1}=\widehat{\mathcal{J}}_{\mathbb{S}_1}.$   
\end{remark}
%%%%%%%%%%%%%%%%%%%%%%%%
In the next result, we present the formula of structured BE when $C={\bf 0}$ and $E={\bf 0}.$
\begin{corollary}\label{coro2:theorem31}
     Suppose $[\x^T,\, \y^T,\z^T]^T$ is  an approximate solution of the  DSPP \eqref{SPP:EQ1} with $A\in \S_n,$ $C={\bf 0},$   $E={\bf 0},$ and $\theta_8, \theta_9, \theta_{10}\neq 0.$  Then, we have 
\begin{align}
    {\bm \eta}^{\mathbb{S}_1}_{\bf sps}(\x,\y,\z)=\left\|\widetilde{\mathcal{J}}_{\mathbb{S}_1}^{T}(\widetilde{\mathcal{J}}_{\mathbb{S}_1}\widetilde{\mathcal{J}}_{\mathbb{S}_1}^T)^{-1}R_{\bm d}\right\|_2,
\end{align}
where $ \widetilde{\mathcal{J}}_{\mathbb{S}_1}\in \R^{(n+m+p)\times \bm{l}}$ is given by
{\footnotesize
\begin{align*}
 \widetilde{\mathcal{J}}_{\mathbb{S}_1}=\bmatrix{ \frac{1}{\theta_1}{\mathcal{K}}_{\x} \Phi_A \mathfrak{D}_{\S_{n}}^{-1}& \frac{1}{\theta_2}{\mathcal{N}}^{n}_{\y}\Phi_B&  {\bf 0} &  -\frac{1}{\theta_8}I_n& {\bf 0} & {\bf 0}\\
 {\bf 0}& \frac{1}{\theta_2}{\mathcal{M}}^m_{\x}\Phi_B & \frac{1}{\theta_5}{\mathcal{N}}^m_{\z}\Phi_D  &  {\bf 0} & -\frac{1}{\theta_9} I_m & {\bf 0} \\
 {\bf 0} & {\bf 0} & \frac{1}{\theta_5}{\mathcal{M}}^{p}_{\y}\Phi_D &  {\bf 0} & {\bf 0}&  -\frac{1}{\theta_{10}}I_p },
\end{align*}}
 $R_f=f-A\x-B^T\y,$ $R_g=g-B\x-D^T\z,$  $R_h=h-D\y,$ and $\bm{l}={\mu}+mn+mp+m+n+p.$
\end{corollary}
\proof Since $C={\bf 0}$ and $E={\bf 0},$ the proof follows by considering $\theta_4=\theta_7=0$. $\square$
%%%%%%%%%%%%%%%%
\begin{remark}
    When $C={\bf 0}$ and $E={\bf 0},$ Lv and Zheng \cite{threeBE2022} derive the structured BE for the  DSPP \eqref{SPP:EQ1}. However, their investigations do not take into account the sparsity pattern of the coefficient matrices. 
\end{remark}
%%%%%%%%%%%%%%%%%%%%%%%%%%%
\section{Derivation of structured BEs for case $(ii)$}\label{sec4}
%%%%%%%%
In this section, we derive explicit formulae for the structured BEs for the  DSPP  for the case $(ii),$ i.e., $A\in \S_n,$ $B_1\neq B_2,$   $C\in \S_m,$ $D_1=D_2=:D$ and $E\in \S_p.$ We use the Lemmas \ref{Lemma1}, \ref{Lemma2} and \ref{sec2:lemma}, and apply a similar methodology used in Section \ref{sec3} to derive the formulae for the structured BEs. In the next result, we present computable formulae for the structured BE ${\bm \eta}^{\mathbb{S}_2}_{\tt sps}(\x,\y,\z)$ by preserving sparsity pattern of the original matrices to the perturbation matrices. Before proceeding, we set the following notations:
\begin{eqnarray*}
\Phi_{B_1}=\diag(\vec(\Theta_{B_1})),~ \Phi_{B_2}=\diag(\vec_{\S}(\Theta_{B_2})),
\end{eqnarray*}
and $\Phi_A, \Phi_C,\Phi_D$ and $\Phi_E$ are same as defined in Section \ref{sec3}.
\begin{theorem}\label{theorem1:sec4}
      Let $[\x^T,\, \y^T,\z^T]^T$ be  an approximate solution of the  DSPP \eqref{SPP:EQ1} with  $A\in \S_n,$ $B_1\neq B_2,$   $C\in \S_m,$ $D_1=D_2=:D,$  $E\in \S_p$ and $ \theta_8, \theta_9,\theta_{10}\neq 0.$  Then, we have 
\begin{align}
    {\bm \eta}^{\mathbb{S}_2}_{\bf sps}(\x,\y,\z)=\left\|\mathcal{J}_{\mathbb{S}_2}^{T}(\mathcal{J}_{\mathbb{S}_2}\mathcal{J}_{\mathbb{S}_2}^T)^{-1}R_{\bm d}\right\|_2,
\end{align}
where $ \mathcal{J}_{\mathbb{S}_2}=[\widetilde{\mathcal{J}}_{\mathbb{S}_2}~~ \mathcal{I}]\in \R^{(n+m+p)\times \bm{l}}$ and $\widetilde{\mathcal{J}}_{\mathbb{S}_2}$ is given by
{\footnotesize
\begin{align*}
 \widetilde{\mathcal{J}}_{\mathbb{S}_2}=\bmatrix{ \frac{1}{\theta_1}{\mathcal{K}}_{\x} \Phi_A\mathfrak{D}_{\S_{n}}^{-1} & \frac{1}{\theta_2}{\mathcal{N}}^{n}_{\y}\Phi_{B_1}&{\bf 0} & {\bf 0} & {\bf 0} & {\bf 0} \\
 {\bf 0}&{\bf 0} & \frac{1}{\theta_3}{\mathcal{M}}^m_{\x}\Phi_{B_2} & -\frac{1}{\theta_4}{\mathcal{K}}_{\y}\Phi_C\mathfrak{D}_{\S_{m}}^{-1} &\frac{1}{\theta_5}{\mathcal{N}}^m_{\z}\Phi_D & {\bf 0} \\
 {\bf 0} & {\bf 0} &{\bf 0} &{\bf 0}& \frac{1}{\theta_5}{\mathcal{M}}^{p}_{\y}\Phi_D & \frac{1}{\theta_7}{\mathcal{K}}_{\z} \Phi_E\mathfrak{D}_{\S_{p}}^{-1}},
\end{align*}}
$R_{\bm d}=[R_f^T,\, R_g^T,\, R_h^T]^T,$ $R_f=f-A\x-B_1^T\y,$ $R_g=g-B_2\x+C\y-D^T\z,$  $R_h=h-D\y-E\z$ and $\bm{l}=\mu+{\sigma}+{\tau}+2mn+mp+m+n+p.$
\end{theorem}    
%%%
\proof Given that $[\widetilde{x}^T,\widetilde{y}^T,\z^T]^T$ is an approximate solution of the  DSPP  \eqref{SPP:EQ1} for the case $(ii).$ Now, we are required to construct perturbations matrices $
       \D A\in \S_n,$ $\D B_1, \D B_2\in \R^{m\times n},$ $
     \D C\in \S_m,$ $\D D\in \R^{p\times m},$   $\D E\in \S_p,$  which maintain the sparsity pattern of $A, B_1, B_2,$ $ C, D, E, $ respectively, and the perturbations $ \D f\in \R^n,$ $ \D g\in \R^m$ and $ \D h\in \R^p$. Using \eqref{ss2:eq24}, $\left(\begin{array}{c}
       \D A,\D B_1,\D B_2,\\
   \D C,  \D D, \D E, \\
     \D f, \D g, \D h
     \end{array}\right)\in\, \mathbb{S}_2$ if and only if $\D A, \D B_1, \D B_2, \D C, \D D, $ $\D E,$ $ \D f, \D g$ and $\D h$ satisfy the following equations:
 \begin{eqnarray}\label{theorem :eq42}
      \left.   \begin{array}{lcl}
       \D A\x+ \D B_1^{T}\y -\D f= R_f, &  \\
       \D B_2\x-\D C\y+\D D^{T}\z- \D g= R_g,& \\
       \D D\y+\D E\z-\D h=R_h,
    \end{array}\right\}
    \end{eqnarray}
    and  $\D A\in \S_n, \D C\in \S_m,$ $\D E\in \S_p.$ \\
     By following a similar the proof methodology  of Theorem \ref{theorem1} and applying Lemma \ref{Lemma1}, we get: \begin{equation}\label{sec4:eq43}
         \mathcal{J}^1_{\mathbb{S}_2}\D X=R_f,~  \mathcal{J}^2_{\mathbb{S}_2}\D X=R_g ~ \mbox{and}~ \mathcal{J}^3_{\mathbb{S}_2}\D X=R_h,
      \end{equation}
      where \begin{eqnarray*}
         && \mathcal{J}^1_{\mathbb{S}_2}=\bmatrix{ \theta_1^{-1}{\mathcal{K}}_{\x}  \Phi_A \mathfrak{D}_{\S_{n}}^{-1}& \theta_2^{-1}{\mathcal{N}}^n_{\y}\Phi_{B_1}&{\bf 0} & {\bf 0} & {\bf 0} & {\bf 0} & -\theta^{-1}_8I_n& {\bf 0} & {\bf 0}}\in \R^{n\times {\bm l}},\\
       && \mathcal{J}^2_{\mathbb{S}_2}= \bmatrix{{\bf 0}& {\bf 0} &\theta_3^{-1}{\mathcal{M}}^m_{\x}\Phi_{B_2} & -\theta_4^{-1}{\mathcal{K}}_{\y}\Phi_C\mathfrak{D}_{\S_{m}}^{-1} &\theta_5^{-1}{\mathcal{N}}^m_{\y}\Phi_D & {\bf 0} &  {\bf 0} & -\theta^{-1}_9 I_m & {\bf 0}}\in \R^{m\times {\bm l}},\\
       &&\mathcal{J}^3_{\mathbb{S}_2}=\bmatrix{{\bf 0} & {\bf 0} &{\bf 0} &{\bf 0}& \theta_5^{-1}{\mathcal{M}}^{p}_{\y}\Phi_D & \theta_7^{-1}{\mathcal{K}}_{\z} \Phi_E\mathfrak{D}_{\S_{p}}^{-1} & {\bf 0} & {\bf 0}&  -\theta^{-1}_{10}I_p }\in \R^{p\times {\bm l}},
      \end{eqnarray*}
      and \begin{eqnarray}
            \D X=\bmatrix{\theta_1\mathfrak{D}_{\S_{n}}\vec_{\S}(\D A\odot \Theta_{A})\\ \theta_2\vec(\D B\odot \Theta_{B_1})\\ 
            \theta_3\vec(\D B\odot \Theta_{B_2})\\\theta_4\mathfrak{D}_{\S_{m}}\vec_{\S}(\D C\odot \Theta_{C})\\ \theta_5\vec(\D D\odot \Theta_D) \\ \theta_7\mathfrak{D}_{\S_{p}}\vec_{\S}(\D E\odot \Theta_{E})\\ \theta_8\D f\\ \theta_9\D g\\ \theta_{10}\D h}\in \R^{\bm l} .
      \end{eqnarray}
      Combining the three equations in \eqref{sec4:eq43}, we obtain \begin{eqnarray}\label{sec4:eq48}
          \mathcal{J}_{\mathbb{S}_2}\D X=R_{\bm d}.
      \end{eqnarray}
      Since,  $\mathcal{J}_{\mathbb{S}_2}$ has full row rank for $\theta_8,\theta_9,\theta_{10}\neq 0$. Therefore, \eqref{sec4:eq48} is consistent and by Lemma \ref{sec2:lemma}, its minimum norm solution is given by 
      \begin{equation}\label{sec4:eq49}
          \D X_{\min}=\mathcal{J}_{\mathbb{S}_2}^T(\mathcal{J}_{\mathbb{S}_2}\mathcal{J}_{\mathbb{S}_2}^T)^{-1}R_{\bm d}.
      \end{equation}
      Now, applying a similar argument to the proof method of Theorem \ref{theorem1}, the required  structured BE is $${\bm \eta}_{\bf sps}^{\mathbb{S}_2}(\x,\y,\z)=\|\D X_{\min}\|_2=\left\|\mathcal{J}_{\mathbb{S}_2}^{T}(\mathcal{J}_{\mathbb{S}_2}\mathcal{J}_{\mathbb{S}_2}^T)^{-1} R_{\bm d}\right\|_2.$$
     Hence, the proof is completed. $\square$
\begin{remark}\label{sec4:remark41}
    The minimal perturbation matrices 
    $\widehat{\D A}_{\tt sps},$  $\widehat{\D C}_{\tt sps},$  $\widehat{\D E}_{\tt sps},$ $\widehat{\D f}_{\tt sps},$ $\widehat{\D g}_{\tt sps}$ and 
    $\widehat{\D h}_{\tt sps}$  can be computed using the formulae provided in Theorem \ref{theorem1} by replacing $\mathcal{J}_{\mathbb{S}_1}$ with $\mathcal{J}_{\mathbb{S}_2}.$ The generating vectors for the minimal perturbation matrices $\widehat{\D B_1}_{\tt sps}$ and
    $\widehat{\D B_2}_{\tt sps}$ are given by 
    $$\vec(\widehat{\D B_1}_{\tt sps})=\frac{1}{\theta_2}\bmatrix{ {\bf 0}_{\mu}& I_{mn}& {\bf 0}_{{\bm{l}}-\mu-mn}}\mathcal{J}_{\mathbb{S}_2}^{T}(\mathcal{J}_{\mathbb{S}_2}\mathcal{J}_{\mathbb{S}_2}^T)^{-1}R_{\bm d} ~~\text{and}$$ 
   $$\vec(\widehat{\D B_2}_{\tt sps})=\frac{1}{\theta_3}\bmatrix{ {\bf 0}_{\mu+mn}& I_{mn}& {\bf 0}_{{\bm l}-\mu-2mn}}\mathcal{J}_{\mathbb{S}_2}^{T}(\mathcal{J}_{\mathbb{S}_2}\mathcal{J}_{\mathbb{S}_2}^T)^{-1}R_{\bm d}.$$
\end{remark}

In the subsequent result, we provide the structured BE while the sparsity pattern is not preserved.
\begin{corollary}\label{SEC4coro:theorem2}
     Let $[\x^T,\, \y^T,\z^T]^T$ be  an approximate solution of the  DSPP \eqref{SPP:EQ1} with $A\in \S_n,$ $C\in \S_m,$   $E\in \S_p,$ and $\theta_8, \theta_9,\theta_{10}\neq 0.$  Then, we have 
\begin{align}
    {\bm \eta}^{\mathbb{S}_2}(\x,\y,\z)=\left\|\widehat{\mathcal{J}}_{\mathbb{S}_2}^{T}(\widehat{\mathcal{J}}_{\mathbb{S}_2}\widehat{\mathcal{J}}_{\mathbb{S}_2}^T)^{-1}R_{\bm d}\right\|_2,
\end{align}
where $ \widehat{\mathcal{J}}_{\mathbb{S}_2}\in \R^{(n+m+p)\times \bm{l}}$ is given by
{\footnotesize
\begin{align*}
 \widehat{\mathcal{J}}_{\mathbb{S}_2}=\bmatrix{ \frac{1}{\theta_1}{\mathcal{M}}_{\x}^n  & \frac{1}{\theta_2}{\mathcal{N}}^{n}_{\y}& {\bf 0} & {\bf 0} & {\bf 0}& {\bf 0} & -\frac{1}{\theta_8}I_n& {\bf 0} & {\bf 0}\\
 {\bf 0}&  {\bf 0}&\frac{1}{\theta_3}{\mathcal{M}}^m_{\x} & -\frac{1}{\theta_4}{\mathcal{K}}_{\y}\mathfrak{D}_{\S_{m}}^{-1} &\frac{1}{\theta_5}{\mathcal{N}}^m_{\z} &  {\bf 0} &  {\bf 0} & -\frac{1}{\theta_9} I_m & {\bf 0} \\
 {\bf 0} & {\bf 0} &{\bf 0}& {\bf 0}& \frac{1}{\theta_5}{\mathcal{M}}^{p}_{\y} & \frac{1}{\theta_7}{\mathcal{K}}_{\z} \mathfrak{D}_{\S_{p}}^{-1} & {\bf 0} & {\bf 0}&  -\frac{1}{\theta_{10}}I_p }.
\end{align*}}
\end{corollary}

%%%%%%%%%%%%%%%%%%%%%%
\proof The proof proceeds by choosing $\Theta_A={\bf 1}_{n\times n},$ $\Theta_{B_1}=\Theta_{B_2}={\bf 1}_{m\times n},$ $\Theta_C={\bf 1}_{m\times m},$ $\Theta_{D}={\bf 1}_{p\times m},$ and  $\Theta_E={\bf 1}_{p\times p}$ in the expression of the structured BE presented in Theorem  \ref{theorem1:sec4}. $\square$
%%%%%%%%%%%%%%%%%%%%%%%%%%%%%
\begin{remark}
    Similar to Corollary \ref{coro2:theorem31}, we can compute the structured BE for the case $(ii)$ with $C={\bf 0}$ and $E={\bf 0}.$ This specific instance of structured BE has also been addressed in \cite{BE2024}. However, our investigation additionally ensures the preservation of the sparsity pattern.
\end{remark}
%%%%%
\section{Derivation of structured BEs for case $(iii)$}\label{sec5}
This section deals with the structured BE of the  DSPP  for the case $(iii),$ i.e., $A\neq A^T$, $C\in \S_m,$ $E\in \S_p,$ $B_1=B_2=:B,$ and $D_1\neq D_2.$ Using a similar technique as in Section \ref{sec3}, in the following theorem, we present the computable formula of the structured BE when sparsity pattern of the original matrices are preserved in the perturbation matrices. Before continuing, we introduce the following notations:
\begin{eqnarray*}
\Phi_{D_1} = \diag(\vec(\Theta_{D_1})), ~
\Phi_{D_2} = \diag(\vec_{\S}(\Theta_{D_2})),
\end{eqnarray*}
along with $\Phi_A$, $\Phi_B$, $\Phi_C$ and $\Phi_E$ as defined in Section \ref{sec3}.
%%%%%%%%%%
\begin{theorem}\label{theorem2}
      Let $[\x^T,\, \y^T,\z^T]^T$ be  an approximate solution of the  DSPP \eqref{SPP:EQ1} with  $C\in \S_m,$   $E\in \S_p,$ and $\theta_8, \theta_9, \theta_{10}\neq 0,$ $B_1=B_2=:B$ and $D_1\neq D_2.$  Then, we have 
\begin{align}
    {\bm \eta}^{\mathbb{S}_3}_{\bf sps}(\x,\y,\z)=\left\|\mathcal{J}_{\mathbb{S}_3}^{T}(\mathcal{J}_{\mathbb{S}_3}\mathcal{J}_{\mathbb{S}_3}^T)^{-1}R_{\bm d}\right\|_2,
\end{align}
$ \mathcal{J}_{\mathbb{S}_3}=[\widetilde{\mathcal{J}}_{\mathbb{S}_3}~~\mathcal{I}]\in \R^{(n+m+p)\times \bm{l}}, $ $\widetilde{\mathcal{J}}_{\mathbb{S}_3}$ is given by

\begin{align*}
\widetilde{\mathcal{J}}_{\mathbb{S}_3}=\bmatrix{ \frac{1}{\theta_1}{\mathcal{M}}_{\x}^n \Phi_A & \frac{1}{\theta_2}{\mathcal{N}}^{n}_{\y}\Phi_B& {\bf 0} & {\bf 0}& {\bf 0}& {\bf 0} \\
 {\bf 0}& \frac{1}{\theta_2}{\mathcal{M}}^m_{\x}\Phi_B & -\frac{1}{\theta_4}{\mathcal{K}}_{\y}\Phi_C\mathfrak{D}_{\S_{m}}^{-1} &\frac{1}{\theta_5}{\mathcal{N}}^m_{\z}\Phi_{D_1} & {\bf 0} &  {\bf 0}  \\
 {\bf 0} & {\bf 0} &{\bf 0}& {\bf 0}& \frac{1}{\theta_6}{\mathcal{M}}^{p}_{\y}\Phi_{D_2} & \frac{1}{\theta_7}{\mathcal{K}}_{\z} \Phi_E\mathfrak{D}_{\S_{p}}^{-1} },
\end{align*}
$R_{\bm d}=[R_f^T,\, R_g^T,\, R_h^T]^T,$ $R_f=f-A\x-B^T\y,$ $R_g=g-B\x+C\y-D_1^T\z,$  $R_h=h-D_2\y-E\z$ and $\bm{l}=n^2+{\sigma}+{\tau}+mn+2mp+m+n+p.$
\end{theorem}
%\vspace{-4mm}
%%%%%%%%%%%%%%%
\proof %The proof follows a similar methodology to that of Theorem \ref{theorem1}. 
For the given approximate solution $[\widetilde{x}^T,\widetilde{y}^T,\z^T]^T,$ we are required to construct perturbations matrices $
       \D A\in \R^n,$ $\D B\in \R^{m\times n},$ $
     \D C\in \S_m,$ $\D D_1,\D D_2\in \R^{p\times m},$   $\D E\in \S_p,$  which maintain the sparsity pattern of $A, B, C, D_1, D_2, E, $ respectively, and the perturbations $ \D f\in \R^n,$ $ \D g\in \R^m$ and $ \D h\in \R^p$. By \eqref{s2:eq24}, $\left(\begin{array}{c}
       \D A,\D B,\D C,\\
     \D D_1, \D D_2, \D E, \\
    \D f,  \D g, \D h
     \end{array}\right)\in\, \mathbb{S}_3$ if and only if $\D A, \D B, \D C, \D D_1, \D D_2$ $\D E,$ $ \D f, \D g$ and $\D h$ satisfy the following equations:
 \begin{eqnarray}\label{theorem2:eq42}
      \left.   \begin{array}{lcl}
       \D A\x+ \D B^{T}\y -\D f= R_f, &  \\
       \D B\x-\D C\y+\D D_1^{T}\z- \D g= R_g,& \\
       \D D_2\y+\D E\z-\D h=R_h,
    \end{array}\right\}
    \end{eqnarray}
    and  $\D C\in \S_m,$ $\D E\in \S_p.$ \\
     By following a similar the proof methodology  of Theorem \ref{theorem1} and applying Lemmas \ref{Lemma1} and \ref{Lemma2}, we get: \begin{equation}\label{eq43}
         \mathcal{J}^1_{\mathbb{S}_3}\D X=R_f,~  \mathcal{J}^2_{\mathbb{S}_3}\D X=R_g ~ \mbox{and}~ \mathcal{J}^3_{\mathbb{S}_3}\D X=R_h,
      \end{equation}
      where \begin{eqnarray*}
         && \mathcal{J}^1_{\mathbb{S}_3}=\bmatrix{ \theta_1^{-1}{\mathcal{M}}_{\x}^n  \Phi_A & \theta_2^{-1}{\mathcal{N}}^n_{\y}\Phi_B& {\bf 0} & {\bf 0} & {\bf 0} & {\bf 0} & -\theta^{-1}_8I_n& {\bf 0} & {\bf 0}}\in \R^{n\times {\bm l}},\\
       && \mathcal{J}^2_{\mathbb{S}_3}= \bmatrix{{\bf 0}& \theta_2^{-1}{\mathcal{M}}^m_{\x}\Phi_B & -\theta_4^{-1}{\mathcal{K}}_{\y}\Phi_C\mathfrak{D}_{\S_{m}}^{-1} &\theta_5^{-1}{\mathcal{N}}^m_{\y}\Phi_{D_1} & {\bf 0} & {\bf 0} &  {\bf 0} & -\theta^{-1}_9 I_m & {\bf 0}}\in \R^{m\times {\bm l}},\\
       &&\mathcal{J}^3_{\mathbb{S}_3}=\bmatrix{{\bf 0} & {\bf 0}&{\bf 0}   &{\bf 0}& \theta_6^{-1}{\mathcal{M}}^{p}_{\y}\Phi_D & \theta_7^{-1}{\mathcal{K}}_{\z} \Phi_E\mathfrak{D}_{\S_{p}}^{-1} & {\bf 0} & {\bf 0}&  -\theta^{-1}_{10}I_p }\in \R^{p\times {\bm l}},
      \end{eqnarray*}
      and \begin{eqnarray}
            \D X=\bmatrix{\theta_1\vec(\D A\odot \Theta_{A})\\ \theta_2\vec(\D B\odot \Theta_B)\\ \theta_4\mathfrak{D}_{\S_{m}}\vec_{\S}(\D C\odot \Theta_{C})\\ \theta_5\vec(\D D\odot \Theta_{D_1}) \\ 
            \theta_6\vec(\D D\odot \Theta_{D_2}) \\\theta_7\mathfrak{D}_{\S_{p}}\vec_{\S}(\D E\odot \Theta_{E})\\ \theta_8\D f\\ \theta_9\D g\\ \theta_{10}\D h}\in \R^{\bm l} .
      \end{eqnarray}
      Combining the three equations in \eqref{eq43}, we obtain \begin{eqnarray}\label{eq48}
          \mathcal{J}_{\mathbb{S}_3}\D X=R_{\bm d}.
      \end{eqnarray}
      Observed that, for $\theta_8,\theta_9,\theta_{10}\neq 0,$ $\mathcal{J}_{\mathbb{S}_3}$ has full row rank. Therefore, \eqref{eq48} is consistent, and by Lemma \ref{sec2:lemma}, its minimum norm solution is given by 
      \begin{equation}\label{eq49}
          \D X_{\min}=\mathcal{J}_{\mathbb{S}_3}^T(\mathcal{J}_{\mathbb{S}_3}\mathcal{J}_{\mathbb{S}_3}^T)^{-1}R_{\bm d}.
      \end{equation}
      Applying a similar  proof technique of Theorem \ref{theorem1}, the required  structured BE is $${\bm \eta}_{\bf sps}^{\mathbb{S}_3}(\x,\y,\z)=\|\D X_{\min}\|_2=\left\|\mathcal{J}_{\mathbb{S}_3}^{T}(\mathcal{J}_{\mathbb{S}_3}\mathcal{J}_{\mathbb{S}_3}^T)^{-1}R_{\bm d}\right\|_2.$$
     Hence, the proof is completed.
      $\square$
%%%%%%%%%%%%%%%%
\begin{remark}\label{remark41}
    The minimal perturbation matrices $\widehat{\D B}_{\tt sps},$ $\widehat{\D C}_{\tt sps},$  $\widehat{\D E}_{\tt sps},$ $\widehat{\D f}_{\tt sps},$ $\widehat{\D g}_{\tt sps}$ and
    $\widehat{\D h}_{\tt sps}$ can be computed using the formulae provided in Theorem \ref{theorem1} with $\mathcal{J}_{\mathbb{S}_1}=\mathcal{J}_{\mathbb{S}_3}.$ The generating vector for the minimal perturbation matrices $\widehat{\D A}_{\tt sps},$ $\widehat{\D D_1}_{\tt sps}$ and $\widehat{\D 
    D_2}_{\tt sps}$    are given by
    \begin{eqnarray*}
       && \vec(\widehat{\D A}_{\tt sps})=\frac{1}{\theta_1}\bmatrix{I_{n^2}& {\bf 0}_{{\bm l}-n^2}}\mathcal{J}_{\mathbb{S}_3}^{T}(\mathcal{J}_{\mathbb{S}_3}\mathcal{J}_{\mathbb{S}_3}^T)^{-1}R_{\bm d}.\\
      && \vec(\widehat{\D D_1}_{\tt sps})=\frac{1}{\theta_5}\bmatrix{ {\bf 0}_{n^2+\m+mn}& I_{mp}& {\bf 0}_{{\tau}+mp+n+m+p}}\mathcal{J}_{\mathbb{S}_3}^{T}(\mathcal{J}_{\mathbb{S}_3}\mathcal{J}_{\mathbb{S}_3}^T)^{-1}R_{\bm d}, ~\text{and}\\
      && \vec(\widehat{\D D_2}_{\tt sps})=\frac{1}{\theta_6}\bmatrix{ {\bf 0}_{n^2+\m+mn+mp}& I_{mp}& {\bf 0}_{\tau+n+m+p}}\mathcal{J}_{\mathbb{S}_3}^{T}(\mathcal{J}_{\mathbb{S}_3}\mathcal{J}_{\mathbb{S}_3}^T)^{-1}R_{\bm d}.
    \end{eqnarray*}
\end{remark}
%%%%%%%%%%%%%%%
In the next result, we present the structured BE when the sparsity pattern is not considered.
\begin{corollary}\label{coro:theorem2}
     For the approximate solution  $[\x^T,\, \y^T,\z^T]^T$ of the  DSPP \eqref{SPP:EQ1} with  $C\in \S_m,$   $E\in \S_p,$ and $ \theta_8, \theta_9, \theta_{10}\neq 0,$   we have 
\begin{align}
    {\bm \eta}^{\mathbb{S}_3}(\x,\y,\z)=\left\|\widehat{\mathcal{J}}_{\mathbb{S}_3}^{T}(\widehat{\mathcal{J}}_{\mathbb{S}_3}\widehat{\mathcal{J}}_{\mathbb{S}_3}^T)^{-1}R_{\bm d}\right\|_2,
\end{align}
where $ \widehat{\mathcal{J}}_{\mathbb{S}_3}\in \R^{(n+m+p)\times \bm{l}}$ is given by
{\footnotesize
\begin{align*}
 \widehat{\mathcal{J}}_{\mathbb{S}_3}=\bmatrix{ \frac{1}{\theta_1}{\mathcal{M}}_{\x}^n  & \frac{1}{\theta_2}{\mathcal{N}}^{n}_{\y}& {\bf 0} & {\bf 0} & {\bf 0}& {\bf 0} & -\frac{1}{\theta_8}I_n& {\bf 0} & {\bf 0}\\
 {\bf 0}& \frac{1}{\theta_2}{\mathcal{M}}^m_{\x} & -\frac{1}{\theta_4}{\mathcal{K}}_{\y}\mathfrak{D}_{\S_{m}}^{-1} &\frac{1}{\theta_5}{\mathcal{N}}^m_{\z} & {\bf 0}& {\bf 0} &  {\bf 0} & -\frac{1}{\theta_9} I_m & {\bf 0} \\
 {\bf 0} & {\bf 0} &{\bf 0}& {\bf 0}& \frac{1}{\theta_6}{\mathcal{M}}^{p}_{\y} & \frac{1}{\theta_7}{\mathcal{K}}_{\z} \mathfrak{D}_{\S_{p}}^{-1} & {\bf 0} & {\bf 0}&  -\frac{1}{\theta_{10}}I_p }.
\end{align*}}
% $R_f=f-A\x-B^T\y,$ $R_g=g-B\x+C\y-D^T\z,$  $R_h=h-D\y-E\z,$ and $\bm{l}=n^2+{\sigma}+{\tau}+mn+mp+m+n+p.$
\end{corollary}
%%%%%%%%%%%%%%%%%%%%%%
\proof The proof follows by taking $\Theta_A={\bf 1}_{n\times n},$ $\Theta_B={\bf 1}_{m\times n},$ $\Theta_C={\bf 1}_{m\times m},$ $\Theta_{D_1}=\Theta_{D_2}={\bf 1}_{p\times m},$ and  $\Theta_E={\bf 1}_{p\times p}$ in the BE expression provided in Theorem \ref{theorem2}. $\square$
\begin{remark}
    Structured BEs for the DSPP \eqref{SPP:EQ1} when $A\in \S_n,$ $C\in\S_m,$ $E\in\S_p,$ $B_1\neq B_2, $  $D_1\neq D_2,$ or, $A\neq A^T,$ $C\in\S_m,$ $E\in\S_p,$ $B_1\neq B_2, $  $D_1\neq D_2$ can be derived in a similar technique used in this section and in Sections \ref{sec3} and \ref{sec4}. As the derivation process is similar, we have not studied them here in detail. 
\end{remark}
%     Theorem \ref{theorem51} gives $\bm{\eta}^{\mathbb{S}_i}(\x,\y,\z)< \, {\bm{ tol}} \Rightarrow \bm{\zeta}^{\mathcal{S}_i}(\x,\y,\z)< \, {\bm{ tol}}$ for 
%  $i=1,2,3,$ and $\bm{\zeta}^{\mathbb{S}_1}(\x,\y,\z) < {\bm{ tol}} \Rightarrow\, \bm{\eta}^{\mathbb{S}_1}(\x,\y,\z)< \sqrt{8}\,{\bm{tol}}$ and  $\bm{\zeta}^{\mathbb{S}_i}(\x,\y,\z) < {\bm{ tol}} \Rightarrow\, \bm{\eta}^{\mathbb{S}_i}(\x,\y,\z)< \sqrt{9}\,{\bm{ tol}}$ for $i=2,3.$  
%%%%%%%%%%%%%%%%%%%%%%%%%%%%
\section{BE for least squares problem with equality constraints (LSE)}\label{sec:LSE}
%%%%%%%%%%%%
The LSE arises in various applications, including the analysis of large-scale structures, signal processing, and solving inequality-constrained least squares problems \cite{LSE1998, LSE2003}. In this section, using our developed framework, we derive the sparsity preserving BE for the LSE problem.
 For given matrices $B\in \R^{m\times n}$ and $D\in \R^{p\times m},$ and the vectors $f\in \R^n$ and $h\in \R^p,$ consider the following LSE problems:
\begin{align}\label{LSE}
    \min_{y}\,\|f-B^Ty\|_2~~ \text{subject to}~~Dy=h.
\end{align}
The LSE problem has a unique solution  if 
\begin{align*}
    \rank(C)=p~~\text{and}~~\rank\left(\bmatrix{B^T\\C}\right)=m.
\end{align*}
To solve the LSE problem \eqref{LSE}, we use the method of Lagrange multipliers. Then, the first-order optimality conditions for the LSE problem lead to the following augmented system:
\begin{align}\label{LSE2}
    \bmatrix{I_n & B^T & \bf 0\\ B &\bf 0& D^T\\ \bf 0& D& \bf 0}\bmatrix{r\\ y\\ \lambda}=\bmatrix{f\\ \bf 0\\ h},
\end{align}
where $\lambda\in \R^p$ is the vector of Lagrange multipliers and $r:=f-B^Ty$.

Since the \((1,1)\) block is an identity matrix, it is not perturbed. Let $[\widetilde{r}^T, \,\widetilde{y}^T, \,\widetilde{\lambda}^T]^T$ be an approximate solution of the system \eqref{LSE2}, i.e., $\widetilde{y}$ is the approximate solution of the LSE problem \eqref{LSE}. Then, we define the BE for the LSE problem as follows:
\begin{align*}%\label{BE:eq21}
       & {\bm \eta}^{LSE}(\widetilde{y})=\displaystyle{\min_{\left(\begin{array}{c}
       \D B,\D D,\\
      \D f, \D h
     \end{array}\right)\in\,\mathcal{G} }}
     \left\|\bmatrix{  \alpha_1\|\D B\|_F& \alpha_2\|\D D\|_F\\   \alpha_3\|\D f\|_2& \alpha_4\|\D h\|_2}\right\|_F,
     \end{align*}
     where
     \begin{align}\label{lse:be2}
&\nonumber\mathcal{G}=\Bigg\{\left(\begin{array}{c}
       \D B,\D D,\\
      \D f, \D h
     \end{array}\right) {\Bigg |} \bmatrix{I_n & (B+\D B)^T &{\bf 0}\\ B+\D B &\bf 0& (D+\D D)^T\\  {\bf 0}& D+\D D & \bf 0}\bmatrix{\widetilde{r}\\ \widetilde{y}\\ \widetilde{\lambda}}=\bmatrix{f+\D f\\ \bf 0\\ h+\D h},\\ 
     &\hspace{6cm}  \D B\in \R^{m\times n}, \D D\in \R^{p\times m}, \D f\in \R^n,\D h\in \R^p\Bigg\}
     \end{align}
     and $\alpha_i,$ $i=1,\ldots, 4,$ are positive real numbers.
     % and $\widetilde{r}=f-B^T\widetilde{y}.$
%%%%%%%%%%%%

When analyzing the BE under sparsity-preserving perturbation matrices on \( B \) and \( D \), we denote it as \( {\bm \eta}^{LSE}_{\bf sps}(\y) \).
\begin{theorem}\label{theorem:LSE}
%      Let $\y$ be an approximate solution of the  LSE problem \eqref{LSE}. Further, let $\mathcal{J}_{\mathcal{G}}\in \R^{(n+m+p)\times \bm{s}},$ where 
% \begin{align*}
% \mathcal{J}_{\mathcal{G}}=\bmatrix{ \frac{1}{\alpha_1}{\mathcal{N}}^{n}_{\y}\Phi_B  & {\bf 0} &-\frac{1}{\alpha_3}I_n & {\bf 0}\\
%  \frac{1}{\alpha_1}M^m_{\widetilde{r}}\Phi_B & \frac{1}{\alpha_2}N^m_{\widetilde{\lambda}}\Phi_D & \bf 0&\bf 0 \\
%   {\bf 0} & \frac{1}{\alpha_2}{\mathcal{M}}^{p}_{\y}\Phi_D & \bf 0&-\frac{1}{\alpha_4}I_p},
% \end{align*}
% $\widehat{R}_{\bm d}=[\widehat{R}_f^T,\, \widehat{R}_g^T,\, \widehat{R}_h^T]^T,$ $\widehat{R}_f=f-\widetilde{r}-B^T\y,$ $\widehat{R}_g=-B\widetilde{r}-D^T\z,$  $\widehat{R}_h=h-D\y,$ and $\bm{s}=mn+mp+n+p.$
% When the condition $\rank(\mathcal{J}_{\mathcal{G}})=\rank([\mathcal{J}_{\mathcal{G}}~~ \widehat{R}_{\bm d}])$  is satisfied  the BE for the LSE problem is given by
% \begin{align}
%     {\bm \eta}^{LSE}_{\bf sps}(\y):=\left\|\mathcal{J}_{\mathcal{G}}^{T}(\mathcal{J}_{\mathcal{G}}\mathcal{J}_{\mathcal{G}}^T)^{-1}\widehat{R}_{\bm d}\right\|_2.
% \end{align}
Let $[\widetilde{r}^T, \,\widetilde{y}^T, \,\widetilde{\lambda}^T]^T$ be an approximate solution of the system \eqref{LSE2}, i.e., $\widetilde{y}$ is the approximate solution of the LSE problem \eqref{LSE}. Define the matrix \( \mathcal{J}_{\mathcal{G}} \in \mathbb{R}^{(n+m+p) \times \bm{s}} \) as  

\[
\mathcal{J}_{\mathcal{G}} =
\bmatrix{ 
\frac{1}{\alpha_1}{\mathcal{N}}^{n}_{\y}\Phi_B  & {\bf 0} &-\frac{1}{\alpha_3}I_n & {\bf 0} \\  
\frac{1}{\alpha_1}M^m_{\widetilde{r}}\Phi_B & \frac{1}{\alpha_2}N^m_{\widetilde{\lambda}}\Phi_D & {\bf 0} & {\bf 0} \\  
{\bf 0} & \frac{1}{\alpha_2}{\mathcal{M}}^{p}_{\y}\Phi_D & {\bf 0} & -\frac{1}{\alpha_4}I_p},
\]  
where \( \bm{s} = mn + mp + n + p \).  
Additionally, define the residual vector \( \widehat{R}_{\bm d} \) as  
\[
\widehat{R}_{\bm d} = \bmatrix{ \widehat{R}_f^T,\, \widehat{R}_g^T,\, \widehat{R}_h^T }^T,  
\]  
where  
\(
\widehat{R}_f = f - \widetilde{r} - B^T\y, \quad  
\widehat{R}_g = -B\widetilde{r} - D^T\z, \quad  
\widehat{R}_h = h - D\y.  
\) 
If the condition  
\[
\rank(\mathcal{J}_{\mathcal{G}}) = \rank\left([\mathcal{J}_{\mathcal{G}} \quad \widehat{R}_{\bm d}]\right)
\]  
is satisfied, then the BE for the LSE problem is given by  

\[
{\bm \eta}^{LSE}_{\bf sps}(\y) = \left\| \mathcal{J}_{\mathcal{G}}^{\dagger}  \widehat{R}_{\bm d} \right\|_2.
\]  
\end{theorem}
\proof Let $[\widetilde{r}^T, \,\widetilde{y}^T, \,\widetilde{\lambda}^T]^T$ is the approximate solution of the augmented system \eqref{LSE2}. Then, we need to find perturbation vectors $\D f\in \R^{n},$ $\D h\in \R^{p},$ and sparsity preserving perturbation matrices $\D B\in \R^{m\times n}$ and $\D D\in \R^{p\times m}$ such that \eqref{lse:be2} holds. Then, we have 
\begin{eqnarray}\label{lse:eq33}
      \left.   \begin{array}{lcl}
        \D B^{T}\y -\D f= \widehat{R}_f, &  \\
       \D B\widetilde{r}+\D D^{T}\widetilde{\lambda}= \widehat{R}_g,& \\
       \D D\y-\D h=\widehat{R}_h.
    \end{array}\right\}
    \end{eqnarray}
    By using a similar method to the proof method of Theorem \ref{theorem1}, we obtain:
    \begin{align}\label{lse66}
        \mathcal{J}_{\mathcal{G}}\D X^{LSE}=\widehat{R}_{\bm d},
    \end{align}
    where $\D X^{LSE}=\bmatrix{\alpha_1\vec(\D B\odot \Theta_B)\\ \alpha_2\vec_{\S}(\D D\odot \Theta_{D})\\ \alpha_3\D f \\\alpha_4\D h}.$
    
   \noindent The linear system \eqref{lse66} is consistent if and only if $\rank(\mathcal{J}_{\mathcal{G}})=\rank([\mathcal{J}_{\mathcal{G}}~~ \widehat{R}_{\bm d}]),$ and the minimum norm solution is given by $\D X^{LSE}_{\min}=\mathcal{J}_{\mathcal{G}}^{\dagger}\widehat{R}_{\bm d}.$ 
  Hence, by applying a similar proof technique to that of Theorem \eqref{theorem1},  the required structured BE is attained.
$\blacksquare$
%%%%%%%%%%%%
\section{BE for least squares problem with equality constraints (LSE)}\label{sec:LSE}
%%%%%%%%%%%%
The LSE arises in various applications, including the analysis of large-scale structures, signal processing, and solving inequality-constrained least squares problems \cite{LSE1998, LSE2003}. In this section, using our developed framework, we derive the sparsity preserving BE for the LSE problem.
 For given matrices $B\in \R^{m\times n}$ and $D\in \R^{p\times m},$ and the vectors $f\in \R^n$ and $h\in \R^p,$ consider the following LSE problems:
\begin{align}\label{LSE}
    \min_{y}\,\|f-B^Ty\|_2~~ \text{subject to}~~Dy=h.
\end{align}
The LSE problem has a unique solution  if 
\begin{align*}
    \rank(C)=p~~\text{and}~~\rank\left(\bmatrix{B^T\\C}\right)=m.
\end{align*}
To solve the LSE problem \eqref{LSE}, we use the method of Lagrange multipliers. Then, the first-order optimality conditions for the LSE problem lead to the following augmented system:
\begin{align}\label{LSE2}
    \bmatrix{I_n & B^T & \bf 0\\ B &\bf 0& D^T\\ \bf 0& D& \bf 0}\bmatrix{r\\ y\\ \lambda}=\bmatrix{f\\ \bf 0\\ h},
\end{align}
where $\lambda\in \R^p$ is the vector of Lagrange multipliers and $r:=f-B^Ty$.

Since the \((1,1)\) block is an identity matrix, it is not perturbed. Let $[\widetilde{r}^T, \,\widetilde{y}^T, \,\widetilde{\lambda}^T]^T$ be an approximate solution of the system \eqref{LSE2}, i.e., $\widetilde{y}$ is the approximate solution of the LSE problem \eqref{LSE}. Then, we define the BE for the LSE problem as follows:
\begin{align*}%\label{BE:eq21}
       & {\bm \eta}^{LSE}(\widetilde{y})=\displaystyle{\min_{\left(\begin{array}{c}
       \D B,\D D,\\
      \D f, \D h
     \end{array}\right)\in\,\mathcal{G} }}
     \left\|\bmatrix{  \alpha_1\|\D B\|_F& \alpha_2\|\D D\|_F\\   \alpha_3\|\D f\|_2& \alpha_4\|\D h\|_2}\right\|_F,
     \end{align*}
     where
     \begin{align}\label{lse:be2}
&\nonumber\mathcal{G}=\Bigg\{\left(\begin{array}{c}
       \D B,\D D,\\
      \D f, \D h
     \end{array}\right) {\Bigg |} \bmatrix{I_n & (B+\D B)^T &{\bf 0}\\ B+\D B &\bf 0& (D+\D D)^T\\  {\bf 0}& D+\D D & \bf 0}\bmatrix{\widetilde{r}\\ \widetilde{y}\\ \widetilde{\lambda}}=\bmatrix{f+\D f\\ \bf 0\\ h+\D h},\\ 
     &\hspace{6cm}  \D B\in \R^{m\times n}, \D D\in \R^{p\times m}, \D f\in \R^n,\D h\in \R^p\Bigg\}
     \end{align}
     and $\alpha_i,$ $i=1,\ldots, 4,$ are positive real numbers.
     % and $\widetilde{r}=f-B^T\widetilde{y}.$
%%%%%%%%%%%%

When analyzing the BE under sparsity-preserving perturbation matrices on \( B \) and \( D \), we denote it as \( {\bm \eta}^{LSE}_{\bf sps}(\y) \).
\begin{theorem}\label{theorem:LSE}
%      Let $\y$ be an approximate solution of the  LSE problem \eqref{LSE}. Further, let $\mathcal{J}_{\mathcal{G}}\in \R^{(n+m+p)\times \bm{s}},$ where 
% \begin{align*}
% \mathcal{J}_{\mathcal{G}}=\bmatrix{ \frac{1}{\alpha_1}{\mathcal{N}}^{n}_{\y}\Phi_B  & {\bf 0} &-\frac{1}{\alpha_3}I_n & {\bf 0}\\
%  \frac{1}{\alpha_1}M^m_{\widetilde{r}}\Phi_B & \frac{1}{\alpha_2}N^m_{\widetilde{\lambda}}\Phi_D & \bf 0&\bf 0 \\
%   {\bf 0} & \frac{1}{\alpha_2}{\mathcal{M}}^{p}_{\y}\Phi_D & \bf 0&-\frac{1}{\alpha_4}I_p},
% \end{align*}
% $\widehat{R}_{\bm d}=[\widehat{R}_f^T,\, \widehat{R}_g^T,\, \widehat{R}_h^T]^T,$ $\widehat{R}_f=f-\widetilde{r}-B^T\y,$ $\widehat{R}_g=-B\widetilde{r}-D^T\z,$  $\widehat{R}_h=h-D\y,$ and $\bm{s}=mn+mp+n+p.$
% When the condition $\rank(\mathcal{J}_{\mathcal{G}})=\rank([\mathcal{J}_{\mathcal{G}}~~ \widehat{R}_{\bm d}])$  is satisfied  the BE for the LSE problem is given by
% \begin{align}
%     {\bm \eta}^{LSE}_{\bf sps}(\y):=\left\|\mathcal{J}_{\mathcal{G}}^{T}(\mathcal{J}_{\mathcal{G}}\mathcal{J}_{\mathcal{G}}^T)^{-1}\widehat{R}_{\bm d}\right\|_2.
% \end{align}
Let $[\widetilde{r}^T, \,\widetilde{y}^T, \,\widetilde{\lambda}^T]^T$ be an approximate solution of the system \eqref{LSE2}, i.e., $\widetilde{y}$ is the approximate solution of the LSE problem \eqref{LSE}. Define the matrix \( \mathcal{J}_{\mathcal{G}} \in \mathbb{R}^{(n+m+p) \times \bm{s}} \) as  

\[
\mathcal{J}_{\mathcal{G}} =
\bmatrix{ 
\frac{1}{\alpha_1}{\mathcal{N}}^{n}_{\y}\Phi_B  & {\bf 0} &-\frac{1}{\alpha_3}I_n & {\bf 0} \\  
\frac{1}{\alpha_1}M^m_{\widetilde{r}}\Phi_B & \frac{1}{\alpha_2}N^m_{\widetilde{\lambda}}\Phi_D & {\bf 0} & {\bf 0} \\  
{\bf 0} & \frac{1}{\alpha_2}{\mathcal{M}}^{p}_{\y}\Phi_D & {\bf 0} & -\frac{1}{\alpha_4}I_p},
\]  
where \( \bm{s} = mn + mp + n + p \).  
Additionally, define the residual vector \( \widehat{R}_{\bm d} \) as  
\[
\widehat{R}_{\bm d} = \bmatrix{ \widehat{R}_f^T,\, \widehat{R}_g^T,\, \widehat{R}_h^T }^T,  
\]  
where  
\(
\widehat{R}_f = f - \widetilde{r} - B^T\y, \quad  
\widehat{R}_g = -B\widetilde{r} - D^T\z, \quad  
\widehat{R}_h = h - D\y.  
\) 
If the condition  
\[
\rank(\mathcal{J}_{\mathcal{G}}) = \rank\left([\mathcal{J}_{\mathcal{G}} \quad \widehat{R}_{\bm d}]\right)
\]  
is satisfied, then the BE for the LSE problem is given by  

\[
{\bm \eta}^{LSE}_{\bf sps}(\y) = \left\| \mathcal{J}_{\mathcal{G}}^{\dagger}  \widehat{R}_{\bm d} \right\|_2.
\]  
\end{theorem}
\proof Let $[\widetilde{r}^T, \,\widetilde{y}^T, \,\widetilde{\lambda}^T]^T$ is the approximate solution of the augmented system \eqref{LSE2}. Then, we need to find perturbation vectors $\D f\in \R^{n},$ $\D h\in \R^{p},$ and sparsity preserving perturbation matrices $\D B\in \R^{m\times n}$ and $\D D\in \R^{p\times m}$ such that \eqref{lse:be2} holds. Then, we have 
\begin{eqnarray}\label{lse:eq33}
      \left.   \begin{array}{lcl}
        \D B^{T}\y -\D f= \widehat{R}_f, &  \\
       \D B\widetilde{r}+\D D^{T}\widetilde{\lambda}= \widehat{R}_g,& \\
       \D D\y-\D h=\widehat{R}_h.
    \end{array}\right\}
    \end{eqnarray}
    By using a similar method to the proof method of Theorem \ref{theorem1}, we obtain:
    \begin{align}\label{lse66}
        \mathcal{J}_{\mathcal{G}}\D X^{LSE}=\widehat{R}_{\bm d},
    \end{align}
    where $\D X^{LSE}=\bmatrix{\alpha_1\vec(\D B\odot \Theta_B)\\ \alpha_2\vec_{\S}(\D D\odot \Theta_{D})\\ \alpha_3\D f \\\alpha_4\D h}.$
    
   \noindent The linear system \eqref{lse66} is consistent if and only if $\rank(\mathcal{J}_{\mathcal{G}})=\rank([\mathcal{J}_{\mathcal{G}}~~ \widehat{R}_{\bm d}]),$ and the minimum norm solution is given by $\D X^{LSE}_{\min}=\mathcal{J}_{\mathcal{G}}^{\dagger}\widehat{R}_{\bm d}.$ 
  Hence, by applying a similar proof technique to that of Theorem \eqref{theorem1},  the required structured BE is attained.
$\square$
%%%%%%%%%%%%
\section{Numerical Experiments}\label{sec:numerical}
%%%%%%%%%%%%%%%%%%
In this section, we carry out several numerical experiments to validate our theoretical findings and the strong backward stability of numerical algorithms for solving the  DSPP.  For Examples \ref{exam2}-\ref{example5}, we consider $\theta_1=\frac{1}{\|A\|_F},$ $\theta_2=\frac{1}{\|B_1\|_F},$ $\theta_3=\frac{1}{\|B_2\|_F},$ $\theta_4=\frac{1}{\|C\|_F},$ $\theta_5=\frac{1}{\|D_1\|_F},$ $\theta_6=\frac{1}{\|D_2\|_F},$ $\theta_7=\frac{1}{\|E\|_F},$ $\theta_8=\frac{1}{\|f\|_2},$   $\theta_9=\frac{1}{\|g\|_2}$ and
  $\theta_{10}=\frac{1}{\|h\|_2}.$ 
   All numerical experiments are performed using MATLAB (version R2023b) on an Intel(R) Core(TM) i7-10700 CPU running at 2.90GHz with 16GB of memory. The machine precision is set to $2.2 \times 10^{-16}$.

\begin{exam}
  We consider the  DSPP \eqref{SPP:EQ1} with the block matrices $A\in \S_5,$  $B\in \R^{3\times 5},$  $C\in \S_3,$ $D\in \R^{2\times 3},$  $E\in \S_2,$ $f\in \R^5,$ $g\in \R^3$ and $h\in \R^2;$ that are given by
  \begin{align*}
     & A=\bmatrix{-0.4083&	0.3472&	0&	0.0636&	0\\
0.3472&	-0.8593&	0.0647&	0.1433&	0\\
0&	0.0647&	0&	0&	0.3129\\
0.0636&	0.1433&	0&	-0.4236&	-1.2123\\
0&	0&	0.3129&	-1.2123&	0}, 
\end{align*}%\\
\begin{align*}
&B=\bmatrix{0&	0&	0&	0&	0.0961\\
-2.2777&	0&	-0.1180&	0&	0\\
1.0582&	0.4363&	0&	1.4115&	-0.0146},~ C=\bmatrix{0&	0&	-0.2299\\
0&	-0.7390&	1.0800\\
-0.2299&	1.0800&	0},\\
&D=\bmatrix{0&	0&	-0.38734\\
0&	0&-0.31964},~ E=\bmatrix{0&	0.5387\\
0.5387&	0},~ f=\bmatrix{-2.5245
\\-1.0063\\
-0.4242\\
-0.6612\\
0.7276
},~ g=\bmatrix{0.4566\\-0.5062\\-1.1846},
%\\
%&\text{and} ~ 
  \end{align*}
 and $h=\bmatrix{0.7818\\-0.0804}.$
 We take the approximate solution $\bm{\widetilde{w}}=[\x^T,\y^T,\z^T]^T,$ where $$\x=\bmatrix{-4.4871\\
11.3517\\
100.3742\\
18.4213\\
-1.6524}, ~\y=\bmatrix{-86.5918\\
5.4127\\
2.6903\\}~ \text{and} ~\z=\bmatrix{1.4512\\
3.3886}$$ with residue $\|\A \widetilde{\bm w}-{\bm d}\|_2=3.4180\times 10^{-02}.$ Further, we take $\theta_i=1,$ for all $i=1,2,\ldots, 10.$ We compute the unstructured BE $\bm{\eta}(\widetilde{\bm w})$ using formulae provided in \eqref{BE:EQ2}, structured BE with preserving sparsity $\bm{\eta}^{\mathbb{S}_1}_{\tt sps}(\x,\y,\z)$ using Theorem \ref{theorem1}, and structured BE without preserving sparsity $\bm{\eta}^{\mathbb{S}_1}(\x,\y,\z)$ using Corollary \ref{coro1:theorem31}. The computed values are given as follows:
\begin{eqnarray*}
  \bm{\eta}(\widetilde{\bm w})= 5.0177\times 10^{-05},~				\bm{\eta}^{\mathbb{S}_1}_{\tt sps}(\x,\y,\z)=2.8084\times 10^{-03}~  \text{and}~\bm{\eta}^{\mathbb{S}_1}(\x,\y,\z)=2.9142\times 10^{-04}.
\end{eqnarray*}
The structure-preserving minimal perturbation matrices that preserve the sparsity pattern as well  are given as follows:
\begin{eqnarray*}
    &\widehat{\D A}_{\tt sps}=10^{-04}\times\bmatrix{-2.1152&	1.4013&0&	4.4380&	0\\
1.4013&	6.4477&	-1.6479&	4.9887	0\\
0&	-1.6479	&0&	0&	10.0133\\
4.4380&	4.9887&	0&	-0.7883&	1.0675\\
0&	0&	10.0133&	1.0675&	0},\\
&\widehat{\D B}_{\tt sps}=10^{-04}\times\bmatrix{0&	0&	0&	0&	10.9193\\
1.1455&	0&	2.6987&	0&	0\\
2.9280&	-2.6709&	0&	-6.9292&	0.9127},
%%%
\end{eqnarray*}
\begin{eqnarray*}
&\widehat{\D C}_{\tt sps}=10^{-05}\times\bmatrix{0&	0&	7.7234\\
0&	-16.9623&	5.7956\\
7.7234&	5.7956&	0},~ \widehat{\D D}_{\tt sps}=10^{-04}\times\bmatrix{0&	0&	-2.3102\\
0&	0&	-6.1506},
\end{eqnarray*}
\begin{align*}
&\widehat{\D E}_{\tt sps}=10^{-04}\times\bmatrix{0&	-2.4376\\
-2.4376&	0},~ \widehat{\D f}_{\tt sps}=10^{-05}\times\bmatrix{-4.7141\\
-5.6800\\
53.1269\\
0.4279\\
-1.1206},~ \\
&\widehat{\D g}_{\tt sps}=10^{-05}\times\bmatrix{124.8013\\
-3.1338\\
3.6990}~ \text{and}~ \widehat{\D h}_{\tt sps}=10^{-04}\times\bmatrix{0.6591\\
1.8203}.
\end{align*}
The above perturbation matrices are computed using the formulae presented in Theorem \ref{theorem1}. Furthermore, the structure-preserving minimal perturbation matrices for which the structured BE $\bm{\eta}^{\mathbb{S}_1}(\x,\y,\z)$ is attained are given by 
\begin{align*}
   & \widehat{\D A}=10^{-06}\times \bmatrix{-4.2896&	5.6947&	48.7625&	8.8353&	-2.2268\\
5.6947&	-1.3596	&-7.9961&	-1.1793&	3.7343\\
48.7625	&-7.9961	&-35.1043	&-3.8939&	32.4331\\
8.8353&	-1.1793	&-3.8939&	-0.2468&	5.9104\\
-2.2268&	3.7343	&32.433&	5.9104&	-1.0584},\\
% \end{align*}
% \begin{align*}
&\widehat{\D B}=10^{-06}\times\bmatrix{-81.1245&	6.1830	&-6.7514&	-5.6364&	-54.8512\\
-3.4239&	21.1045&	190.4498&	35.2275&	0.3003\\
6.7250&	-10.8292&	-93.8457&	-17.0866&	3.2526}, \\
&\widehat{\D C}=10^{-06}\times\bmatrix{-31.9502&	83.9647	&-39.5777\\
83.9647	&-10.3722&	-0.0727\\
-39.5777&	-0.0727&	2.4901},\\
&\widehat{\D D}=10^{-06}\times\bmatrix{17.4919&	1.6540&	-1.9033\\
25.4210&	4.8262&	-3.9650},\\
& \widehat{\D E}=10^{-06}\times\bmatrix{-0.3021&	-0.5762\\
-0.5762&-1.0437}, ~ \widehat{\D f}=10^{-07}\times\bmatrix{-9.5598\\
1.1977\\
3.4973\\
0.1340\\
-6.4049},\\
%%%%%%%%%
&\widehat{\D g}=10^{-07}\times\bmatrix{3.6897\\
-19.1626\\
9.2558}~\text{and}~ \widehat{\D h}=10^{-07}\times\bmatrix{2.0819\\
3.0801}.
\end{align*}
 \end{exam}
 %%%%%
 \begin{exam}\label{exam2}
     To test the strong backward stability of numerical algorithms,  we consider the  DSPP \eqref{SPP:EQ1} with 
     \begin{align*}
         &A=GPG(1:3,1:3)\in \S_3,~ B=D=\bmatrix{0&0&1\\0&1&0\\10^4 &0 &0}\in\R^{3\times 3},~ C=\bmatrix{1&-2 &1\\ -2& 6& 0\\ 1& 0&0}\in\R^{3\times 3}, ~\\
         &E=GPG(4:6,4:6)\in \S_3,~ f=\bmatrix{10^8\\10\\0}\in \R^3~\text{and}~ g=h=\bmatrix{10^{-8} \\0\\0}\in \R^3,
     \end{align*}
     where $G=10^{6}\times\diag(1,5,10,50,100,500)$ and $P=[p_{ij}]\in \R^{6\times 6},$ $p_{ij}=\dm{\frac{(i+j-1)!}{(i-1)!(j-1)!}}.$
     The approximate solution $\widetilde{\bm w}=[\x^T,\y^T,\z^T]^T$ of this DSPP is obtained using Gaussian elimination with partial pivoting (GEP), where 
     \begin{eqnarray*}
         \x=10^{-05}\times\bmatrix{60.0120\\-8.0016\\1.0002},~ \y=\bmatrix{6.0012\\2.0004\\-2.0004}~ \text{and}~\z=10^{-13}\times\bmatrix{-1.7109\\0.8556\\-0.0475}.
     \end{eqnarray*}
     We compute the unstructured BE ${\bm \eta}(\widetilde{\bm w})$, structured BEs $ {\bm \eta}_{\tt sps}^{\mathbb{S}_1}(\x,\y,\z)$ and ${\bm \eta}^{\mathbb{S}_1}(\x,\y,\z)$ using the formulae given in \eqref{BE:EQ2}, Theorem \ref{theorem1} and Corollary \ref{coro1:theorem31}, respectively. The obtained BEs are given by 
     \begin{equation}\label{example2:eq1}
         {\bm \eta}(\widetilde{\bm w})=6.9314\times 10^{-27},~ {\bm \eta}_{\tt sps}^{\mathbb{S}_1}(\x,\y,\z)=5.4649\times10^{-06},~ {\bm \eta}^{\mathbb{S}_1}(\x,\y,\z)=4.7907\times 10^{-06}.
     \end{equation}
     From \eqref{example2:eq1}, we can observe that ${\bm \eta}(\widetilde{\bm w})$ of $\mathcal{O}(10^{-27})$ indicates that GEP is backward stable for solving this DSPP. On the other side ${\bm \eta}_{\tt sps}^{\mathbb{S}_1}(\x,\y,\z)$ and ${\bm \eta}^{\mathbb{S}_1}(\x,\y,\z)$ is much bigger than ${\bm \eta}(\widetilde{\bm w})$ implies that GEP for solving this  DSPP is not strongly backward stable. 
     %This shows that a backward stable algorithm for solving this  DSPP is not strongly backward stable.
 \end{exam}
 %%%%%%%%%%%%%%%
 \begin{exam}\label{exam3}
     In this example, we perform a comparison among our obtained structured BEs and the structured BE considered in \cite{BE2024}. For this, we consider the DSPP \eqref{SPP:EQ1} with
     \begin{align*}
         &A=\bmatrix{0.0968&	0&	-0.2438&	-0.2823\\
0& 0&	1.1180&	-1.1611\\
-0.2438&	1.1180&	1.6014&	-0.8693\\
-0.2823&	-1.1611&	-0.8693	&-0.4914},~ B_1=\bmatrix{0&	0&	0.7090&	0\\
1.9046&	0.0928&	-0.0430&	0.0508},~\\
&B_2=\bmatrix{-0.2592&	0&	0.2543&	0.1248\\
0.0876&	1.1375&	0&	0.0766}, ~ C={\bf 0}_{2\times2}, ~D=\bmatrix{0&	1.8070\\
1.0365&	-1.5516} ~ \text{and}~E={\bf 0}_{2\times2}.
     \end{align*}
 Here, $n=4,m=2$ and $p=2.$    Further, we consider the right-hand side vector ${\bm d}=[f^T,g^T,h^T]^T\in \R^8,$ where 
     \begin{align*}
         f=\bmatrix{-1.1251\\
-1.9000\\
-0.4320\\
-1.1422},~ g=\bmatrix{-0.5516\\
1.8738}~ \text{and}~ h=\bmatrix{0.4982\\
0.8347}.
     \end{align*}
     The computed solution using the the MATLAB `\textit{blackshash}' command is $\widetilde{\bm w}=[\x^T,\y^T,\z^T]^T,$ where
      \begin{align*}
         \x=\bmatrix{-1.6927\\
-1.5778\\
1.9746\\
3.5598},~ g=\bmatrix{1.2180\\
0.2757}~ \text{and}~ h=\bmatrix{0.3571\\
-1.8683}.
     \end{align*}
     The computed solution $\widetilde{\bm w}$ has residue  $\|\A\widetilde{\bm w}-{\bm d}\|=1.4864\times 10^{-15}.$  The unstructured BE computed using the formula \eqref{BE:EQ2} is $5.2700\times 10^{-17}$, structured BE using the Theorem $3.2$ of \cite{BE2024} is $4.1137\times 10^{-16},$ structured BE with sparsity using Theorem \ref{theorem2} is $2.7992\times 10^{-16}$ and the structured BE without sparsity using Corollary \ref{coro:theorem2} is $2.5525\times 10^{-16}.$ We observe that all the computed BEs are in unit round-off error and the structured BEs are only one order larger than the unstructured ones. Furthermore, the structured BEs derived in our work and those obtained in the reference \cite{BE2024} exhibit uniform order.  This shows the reliability of our derived structured BEs formulae. One notable advantage of our derived formulae lies in our ability to preserve the sparsity pattern within the perturbation matrices.
 \end{exam}
 %%%%%%%%%%%%%%%
 \begin{exam}\label{example4:GMRES}
     To test the strong backward stability of the GMRES method, in this example, we consider the DSPP \eqref{SPP:EQ1} \cite{HuangNA} with the block matrices 
		\begin{eqnarray}
			\noin A= \bmatrix{I\otimes Z+Z\otimes I &\bf 0\\ \bf 0&I\otimes Z+ Z\otimes I}\in \R^{2r^2\times2r^2}, \,\, B=\bmatrix{I\otimes H& H\otimes I}\in \R^{r^2\times 2r^2}, 
		\end{eqnarray}
	$  D= G\otimes H\in \R^{r^2\times r^2}$ and $C=E={\bf 0}_{r^2\times r^2},$  where $Z=\frac{1}{(r+1)^2}\, \mathrm{tridiag}(-1,2,-1)\in \R^{r\times r},\quad H=\frac{1}{r+1}\,  \mathrm{tridiag}(0,1,-1)\in \R^{r\times r}$  and  $ G=\diag(1, r+1, \ldots, r^2-r+ 1)\in \R^{r\times r}.$ Here, $G\otimes H$ represents the Kronecker product of matrices $G$ and $H$ and $ \mathrm{tridiag}(d_1,d_2,d_3)\in \R^{r\times r}$ represents the 
 $r \times r$ tridiagonal matrix with the subdiagonal entry $d_1$, diagonal entry $d_2$,  and superdiagonal entry $d_3.$  For this problem, the dimension of  $\A$ is $4r^2.$ We use GMRES method to solve this DSPP with termination criteria ${\frac{\|\A {\bm w}_k-{\bm d}\|_2}{\|{\bm d}\|_2}}< {\bm{ tol}},$ where ${\bm w}_k$ is solution at each iterate and  ${\bm{ tol}}=10^{-13}$ and the initial guess vector zero. We compute the structured and unstructured BEs for the solution at the final iteration. The computed BEs for different values of $r$ are listed in Table \ref{tab1}.
     \begin{table}[ht!]
		\centering
		\caption{Values of structured and unstructured BEs of the approximate solution obtained using GMRES for Example \ref{example4:GMRES}.}\label{tab1}
		\resizebox{11cm}{!}{
			\begin{tabular}{@{}ccccc@{}}
				\toprule
				$r$ & ${\frac{\|\A {\bm w}_k-{\bm d}\|_2}{\|{\bm d}\|_2}}$ & ${\bm \eta}(\widetilde{\bm w})$ &  ${\bm \eta}^{\mathbb{S}_1}(\x,\y,\z)$ & ${\bm \eta}^{\mathbb{S}_1}_{\tt sps}(\x,\y,\z)$\\ [1ex]
				\midrule
				%$2$&$9.6852$e-$16$&	$1.2467$e-$16$	&$4.0234$e-$16$&	$6.5433$e-$16$\\
$4$&$1.0593$e-$15$&	$4.1823$e-$17$	&$1.3757$e-$16$	&$4.9831$e-$16$\\
$6$&$2.4960$e-$14$&	$5.3825$e-$16$&	$1.8436$e-$15$&	$7.3871$e-$15$\\
$8$&$2.0868$e-$14$&	$3.0086$e-$16$&	$9.6476$e-$16$&	$5.4875$e-$15$\\
$10$&$3.2981$e-$14$&	$3.4862$e-$16$&	$1.3781$e-$15$&	$9.1775$e-$15$
				 \\ [1ex]
				\toprule
		\end{tabular}}
	\end{table}
 
 From Table \ref{tab1}, we observe that unstructured BE ${\bm \eta}(\widetilde{\bm w}),$ structured BE with preserving sparsity ${\bm \eta}^{\mathbb{S}_1}_{\tt sps}(\x,\y,\z)$ and structured BE without sparsity ${\bm \eta}^{\mathbb{S}_1}(\x,\y,\z)$ are all around order of unit round-off error. Using our obtained structured BEs, we successfully demonstrate that the GMRES method for solving this DSPP exhibits strong backward stability.
 \end{exam}
 %%%%%%%%%
 \begin{exam}\label{example5}
 To test the strong backward stability of the GMRES method to solve the  DSPP for case $(iii)$, in this example, we consider the data matrices for \eqref{SPP:EQ1} as follows:
\begin{eqnarray*}
    &A=randn(n,n), ~ B=sprandn (m,n,0.5),~ 
C=0.5(C_1+C_1^T), ~ \\
&D1=sprandn(p,m,0.5),~ D2=sprandn(p,m,0.5),~E=0.5(E_1+E_1^T), 
\end{eqnarray*}
where $C_1=sprandn(m,m,0.2),$ $E_1=sprandn(p,p,0.3),$ $n=2k,$ $m=k,$ and $p=k.$ Furthermore, we take ${\bm d}=[f^T, g^T, h^T]^T$  such that the exact solution of the DSPP  is $[1,1,\ldots,1]^T\in \R^{n+m+p}.$ Here, the symbols `$randn(n,m)$' and `$sprandn(n,m,\omega)$' stand for normally distributed random matrix and the sparse random matrix with density $\omega$, respectively, of size $n\times m$.

     \begin{figure}
         \centering
        \includegraphics[width=0.5\textwidth]{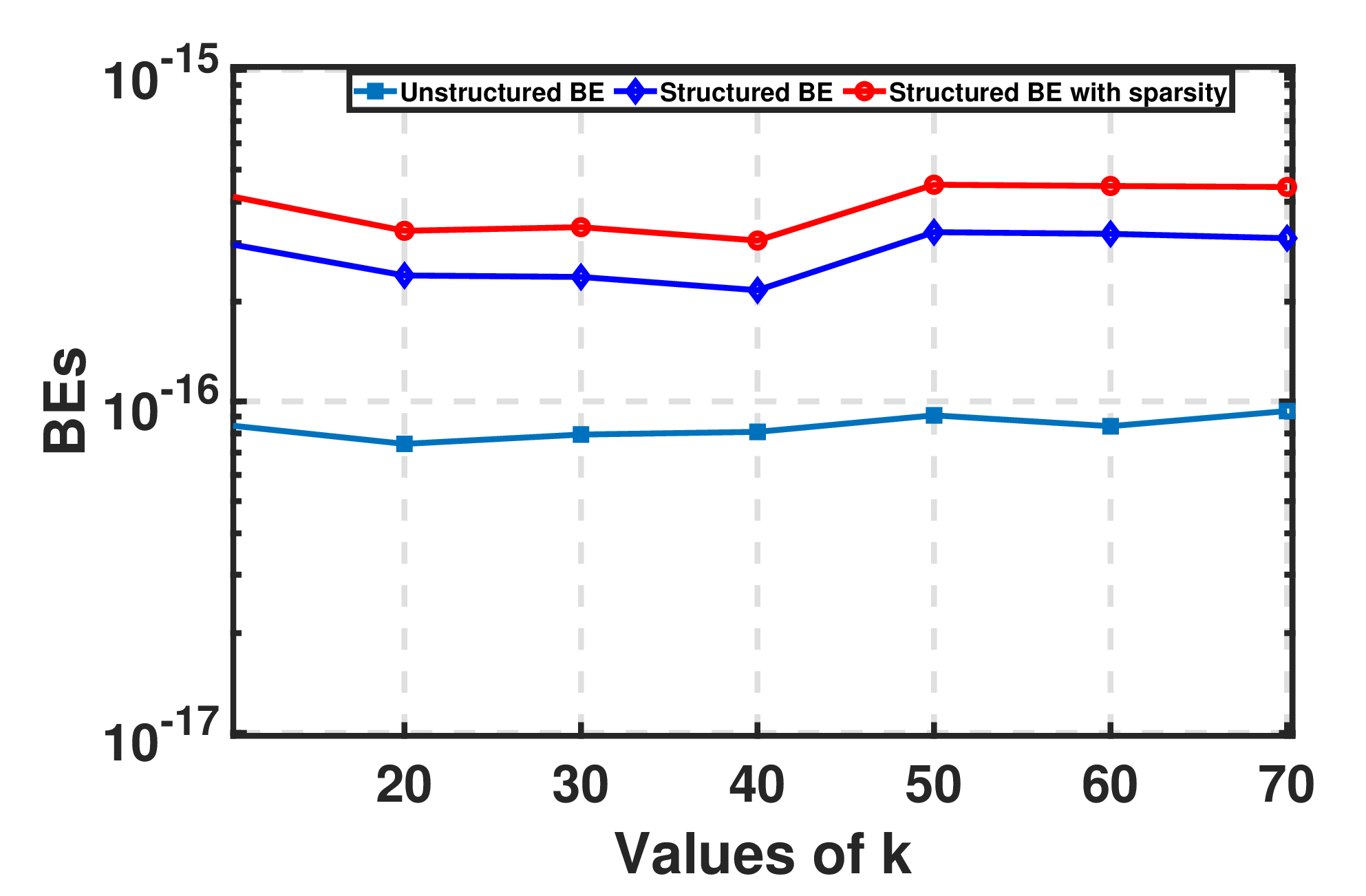}
%{Figures/Strong_example2.eps}
         \caption{Unstructured and structured BEs versus $k$ for Example \ref{example5}.}
         \label{fig:strong}
     \end{figure}

     We employ the GMRES method \cite{gmres} with an initial guess vector set to zero, termination criteria ${\frac{\|\A {\bm w}_k-{\bm d}\|_2}{\|{\bm d}\|_2}}< {\bm{ tol}},$ where ${\bm w}_k$ is solution at each iterate and the tolerance is $\bm{tol}=10^{-13}$. We take values of $k$ from $10$ to $70$ with step size $10.$  For these values of $k$, we plot the unstructured BE $\bm{\eta}(\widetilde{{\bm w}})$ using the formula \eqref{BE:EQ2} (abbreviated as  `\textbf{unstructured BE}'), structured BE $\bm{\eta}_{\tt sps}^{\mathbb{S}_3}(\x, \y,\z)$ (abbreviated as `\textbf{structured BE with sparsity}') using Theorem \ref{theorem2},  and $\bm{\eta}^{\mathbb{S}_3}(\x, \y,\z)$ (abbreviated as `\textbf{structured BE}') using Corollary \ref{coro:theorem2}. From Figure \ref{fig:strong}, it is seen that, for all values of $k,$ the unstructured BE $\bm{\eta}(\widetilde{\bm w})$ are of order $\mathcal{O}(10^{-16})$  and the structured BEs  $ \bm{\eta}_{\tt sps}^{\mathbb{S}_3}(\x, \y,z) $ and $\bm{\eta}^{\mathbb{S}_3}(\x, \y,z)$ are of $\mathcal{O}(10^{-15})$, which are very small.  Notably,  both the structured BEs are only one order less than the unstructured BE.  
     Thus, the utilization of our structured BE affirms that the GMRES method for solving this DSPP is strongly backward stable.
     % This illustrates that the $GMRES$ method for solving the  DSPP is both backward stable and strongly backward stable. 
     Furthermore, the resultant approximate solution corresponds to an exact solution of a nearly perturbed  DSPP of the form \eqref{SPP:EQ1}, which preserves the inherent matrix structure and sparsity pattern.
 \end{exam}
\section{Conclusions}\label{sec:conclusion}
This paper investigates the structured BEs for a class of  DSPP by preserving the inherent matrix structure and sparsity pattern in the perturbation matrices. We derive explicit formulae for the structured BEs (in three cases) and present concrete formulae for structure-preserving minimal perturbation matrices. %for which an approximate solution become an exact solution of a  DSPP.
These perturbation matrices yield a nearly perturbed  DSPP. Thereby, the approximate solution becomes its exact solution. Moreover, our derived framework is used to derive BE for the LSE problem. Several numerical experiments are performed to validate our obtained theory and to test the strong backward stability of numerical algorithms. Our observations reveal that numerical algorithms (such as GEP) demonstrating backward stability do not always exhibit strong backward stability.
%%%%%%%%%%%%%%%%%%%%%%%%%%
\section*{Acknowledgments}
Pinki Khatun expresses her sincere gratitude to the Council of Scientific and Industrial Research (CSIR), New Delhi, India, for the financial support provided through a Research Fellowship (File No. 09/1022(0098)/2020-EMR-I). This work was carried out while she was a Ph.D. scholar at the Indian Institute of Technology Indore, India.
 %%%%%%%%%%%%%%%%%%%%%%%%%%%%%%%
 
\bibliographystyle{abbrvnat}
\bibliography{reference}

 \end{document}